\documentclass[twoside,12pt]{article}
\usepackage{amsmath}
\usepackage[usenames,dvipsnames]{color}

\usepackage{cite}
\usepackage{amsmath}
\usepackage{amsthm}
\usepackage{amssymb}

\definecolor{rosso}{rgb}{0.8,0,0}

\def\gianni #1{{\color{red}#1}}
\def\gian #1{{\color{magenta}#1}}
\def\juerg #1{{\color{Green}#1}}
\def\pier #1{{\color{rosso} #1}}

\def\gianni #1{#1}
\def\gian #1{#1}
\def\juerg #1{#1}
\def\pier #1{#1}

\newcommand{\dn}{\partial_{\bf n}}
\newcommand{\pt}{\partial_t}
\newcommand{\beq}{\begin{equation}}
\newcommand{\eeq}{\end{equation}}
\newcommand{\beqa}{\begin{eqnarray}}
\newcommand{\eeqa}{\end{eqnarray}}

\newcommand{\oma}{\Omega}

\newcommand{\bu}{\overline{u}}
\newcommand{\bmu}{\overline{\mu}}
\newcommand{\br}{\overline{\rho}}
\newcommand{\blam}{\overline{\lambda}}
\newcommand{\uad}{{\cal U}_{\rm ad}}
\newcommand{\lzo}{L^2(\Omega)}
\newcommand{\lio}{L^\infty(\Omega)}

\newcommand{\lzq}{L^2(Q)}
\newcommand{\liq}{L^\infty(Q)}
\newcommand{\ioma}{\int_\Omega}
\newcommand{\txinto}{\int_0^t\!\!\!\int_\Omega}
\newcommand{\texinto}{\int_0^T\!\!\!\int_\Omega}
\newcommand{\tint}{\int_0^t}

\newcommand{\ula}{u^\lambda}
\newcommand{\rla}{\rho^\lambda}
\newcommand{\mula}{\mu^\lambda}
\newcommand{\yla}{y^\lambda}
\newcommand{\zla}{z^\lambda}

\newcommand{\cs}{{\cal S}}
\newcommand{\cx}{{\cal X}}
\newcommand{\cy}{{\cal Y}}

\newcommand{\lzht}{L^2(0,t;H)}
\newcommand{\lzvt}{L^2(0,t;V)}

\newcommand{\dx}{\,dx}
\newcommand{\dy}{\,dy}
\newcommand{\dt}{\,dt}
\newcommand{\ds}{\,ds}

\def\aand{\quad\hbox{and}\quad}
\def\lhs{left-hand side}
\def\rhs{right-hand side}
\def\sfw{straightforward}
\def\accorpa #1#2{\eqref{#1}--\eqref{#2}}
\def\smfrac{\textstyle\frac}

\def\Beq{\begin{equation}}
\def\Eeq{\end{equation}}
\def\Bsist{\begin{eqnarray}}
\def\Esist{\end{eqnarray}}
\let\non\nonumber
\def\step #1 {\vspace{3mm}\underline{\sc Step #1:}  \,\,\,}
\def\Aeps{{\cal A}^\eps}
\def\calQ{{\cal Q}}
\def\Reps{{\cal R}^\eps}
\def\calV{{\cal V}}
\def\calH{{\cal H}}
\def\calL{{\cal L}}

\let\eps\varepsilon
\let\phi\varphi
\def\mueps{\mu^\eps}
\def\rhoeps{\rho^\eps}
\def\xieps{\xi^\eps}
\def\etaeps{\eta^\eps}
\def\aeps{a^\eps}
\def\lambdaeps{\lambda^\eps}
\def\peps{p^\eps}
\def\qeps{q^\eps}
\def\ceps{c_\eps}
\def\bdelta{b^\delta}
\def\feps{f^\eps}
\def\fed{f^{\eps\delta}}
\def\xied{\xi^{\eps\delta}}
\def\zed{\zeta^{\eps\delta}}
\def\fed{f^{\eps\delta}}
\def\Adelta{A_\delta}
\def\phieps{\phi_\eps}
\def\psieps{\psi_\eps}
\def\Meps{M_\eps}

\let\intQt\txinto
\let\intQ\texinto
\def\bintQt{\int_t^T\!\!\!\!\iO}
\def\iO{\int_\Omega}
\def\iot{\int_0^t}
\def\ioT{\int_0^T}
\def\itt{\int_t^T}
\def\<#1>{\mathopen\langle #1\mathclose\rangle}
\def\norma #1{\mathopen \| #1\mathclose \|}
\def\normaV #1{\norma{#1}_V}
\def\normaH #1{\norma{#1}_H}

\def\div{\mathop{\rm div}\nolimits}

\def\genspazio #1#2#3#4#5{#1^{#2}(#5,#4;#3)}
\def\spazio #1#2#3{\genspazio {#1}{#2}{#3}T0}

\def\L {\spazio L}
\def\H {\spazio H}

\def\Lx #1{L^{#1}(\Omega)}

\def\LQ #1{L^{#1}(Q)}

\def\Lsei{\Lx6}

\def\LQ #1{L^{#1}(Q)}

\newcommand{\rz}{\mathbb{R}}
\newcommand{\nz}{\mathbb{N}}
\let\erre\rz

\begin{document}

\begin{center}
{\bf {\huge Distributed optimal control\\[1mm]  of a nonstandard nonlocal \\[5mm] 
phase field system\footnote{This work received \pier{a partial support from the GNAMPA} (\pier{Gruppo} Nazionale per l'Analisi
Matematica, la Probabilit\`{a} e loro Applicazioni) of INDAM (Istituto Nazionale
di Alta \pier{Matematica}) and the \pier{IMATI -- C.N.R. Pavia} for PC and GG.}}}

\vspace{9mm}
{\large Pierluigi Colli$^{\!\dagger}$,
Gianni Gilardi\footnote{Dipartimento di Matematica  ``F. Casorati'',
Universit\`a di Pavia, \pier{Via Ferrata 5,} 27100 Pavia, Italy
(e-mail: pierluigi.colli@unipv.it, gianni.gilardi@unipv.it) },\\[1mm]
and J\"urgen Sprekels\footnote{Weierstrass Institute for 
Applied Analysis and Stochastics,
Mohrenstra\ss e 39, 10117 Berlin and Department of Mathematics, 
Humboldt-Universit\"at zu Berlin, Unter den Linden 6, 10099 Berlin, \pier{Germany}
(e-mail: juergen.sprekels@wias-berlin.de) }}\\[6mm]{\small {\bf Key words:} Distributed optimal control, nonlinear phase field systems,\\ nonlocal operators, 
first-order necessary optimality conditions.\\[2mm]
{\bf AMS (MOS) Subject Classification:} 35K55, 49K20, \pier{74A15.}}
\end{center}

\vspace{3mm}
\begin{abstract} \noindent
We investigate a distributed optimal control problem for a nonlocal phase field 
model of \pier{viscous} Cahn--Hilliard type. The model constitutes a 
nonlocal version of a model for two-species phase segregation  
on an atomic lattice under the presence of diffusion \pier{that has been
studied in a series of papers by P. Podio-Guidugli and the present authors}. The model consists of a highly nonlinear parabolic equation
coupled to an ordinary differential equation. The latter equation contains both 
nonlocal and singular terms that render the analysis difficult. Standard arguments of 
optimal control theory do not apply directly, although the control constraints and the cost functional are of standard type. 
We show that the  problem admits a solution, and we derive the first-order
necessary conditions of optimality. \end{abstract}


\thispagestyle{empty}
\pagestyle{myheadings}
\newcommand\testopari{\sc \pier{Colli \ --- \ Gilardi \ --- \ Sprekels}}
\newcommand\testodispari{\sc \pier{Optimal control of a phase field system}}
\markboth{\testodispari}{\testopari}



\section{Introduction}
Let $\Omega\subset\rz^3$ denote an open and bounded domain whose smooth
boundary $\Gamma$ has the outward unit normal ${\bf n}$\pier{;} let $T>0$ be a given final time, and 
\pier{set} \gianni{$Q:=\oma\times (0,T)$} and
$\Sigma:=\Gamma\times (0,T)$. We study in this paper distributed 
optimal control problems of the following form:

\vspace{3mm}
{\bf (CP)} \,\,Minimize the cost functional
\begin{eqnarray}
\label{cost}
 J(u, \rho,\mu)&\!\!=\!\!&\frac {\beta_1} 2 \texinto |\rho-\rho_Q|^2\dx\dt
\,+\, \frac{\beta_2}{2}\texinto |\mu-\mu_Q|^2 \,dx\,dt\nonumber\\  
&&+\,\frac{\beta_3}{2}\int_0^T\!\!\!\int_\Omega |u|^2\,dx\,dt
\end{eqnarray}
subject to the state system
\begin{eqnarray}
\label{ss1}
(1+2\,g(\rho))\,\pt\mu+\mu\,g'(\rho)\,\pt\rho-\Delta\mu=u\,\quad\mbox{\pier{a.\,e.} in }\,Q,\\
\label{ss2}
\pt\rho\,+\,B[\rho]\,+\,F'(\rho)=\mu\,g'(\rho)\,\quad\mbox{a.\,e. in }\,Q,\\
\label{ss3}
\dn  \mu=0\,\quad\mbox{a.\,e. on }\Sigma,\\
\label{ss4}
\rho(\cdot,0)=\rho_0\,,\quad \mu(\cdot,0)=\mu_0,\,\quad\mbox{a.\,e. in } \,\Omega,
\end{eqnarray}
and to the control constraints
\begin{align}
  \label{uad}
u\in\uad:=\bigl\{u\in H^1(0,T;\lzo)\,:\, \ &0\le u\le u_{\rm max} 
\,\,\,\,\mbox{a.\,e. in }\,Q  \nonumber\\
&\quad \pier{\mbox{and }\,\|u\|_{H^1(0,T;\lzo)}\,\le\,R \bigr\}.}
\end{align}

Here, \gianni{$\beta_1, \beta_2,\beta_3\geq0$ and $R>0$ are given constants\pier{, 
with  $\beta_1 + \beta_2 + \beta_3 >0 $,}         and 
the threshold function $u_{\rm max}\in\liq$ is nonnegative}. 
Moreover, $\rho_Q,\mu_Q\in L^2(Q)$ represent prescribed target functions of the tracking-type functional~$J$. 
Although  more general cost functionals could be admitted 
for large parts of the subsequent analysis, we restrict ourselves to the above situation
for the sake of a simpler exposition.

The state system (\ref{ss1})--(\ref{ss4}) constitutes a {\em nonlocal} version of a
phase field model of Cahn--Hilliard type describing  phase segregation of 
two species (atoms and vacancies, say) on a lattice, which was recently studied in
\cite{CGS3}. In the (simpler) original {\em local} model, which was introduced in \cite{PG},
the nonlocal term $B[\rho]$ is a replaced by the diffusive term $\,-\Delta\rho$.
The local  model has been the subject of intensive research in the past years; in this connection, we refer the reader to \pier{\cite{CGKPS,CGKS1,CGKS2,CGPS3,CGPS6,CGPS7, CGPS4,CGPS5}}. 
In particular, in \cite{CGPSco} the analogue of the control
problem {\bf (CP)} for the local case was investigated for the special \gianni{situation} $g(\rho)=\rho$, 
for which \gianni{the} optimal boundary control problems was studied in~\cite{CGSco1}. 

The state variables of the model are the {\em order parameter} $\rho$, interpreted as a
volumetric density, and the \emph{chemical potential} $\mu$; for physical reasons,
we must have $0\le \rho\le 1$ and $\mu > 0$ almost everywhere in $Q$. The control
function $u$ on the right-hand side of (\ref{ss1}) plays the role of a 
{\em microenergy source}. We remark at this place that the requirement encoded in the definition
of~$\uad$, namely that $u$ be nonnegative, is indispensable for the analysis of the forthcoming sections. Indeed,
it is needed to guarantee the nonnegativity of the chemical potential $\mu$.
 
The nonlinearity $F$ is a double-well potential defined in \pier{the interval} $(0,1)$, whose derivative 
$F'$ is singular at the endpoints $\rho=0$ and $\rho=1$: e.\,g., $F=F_1+F_2$, where
$F_2$ is smooth and 
\beq\label{logpot}
F_1(\rho)=\hat c\,(\rho\,\log(\rho)+(1-\rho)\,\log(1-\rho)), \,\mbox{ with 
a constant $\,\hat c>0$}.
\eeq
The  presence of the nonlocal term $B[\rho]$ in \eqref{ss2} constitutes the main 
difference to the local model. Simple
examples are given by \gian{integral} operators of the form
\beq\label{intop1}
B[\rho](x,t)=\gian{\tint\iO k(t,s,x,y)\,\rho(y,s)\,ds\,dy}
\eeq
and \gian{purely} spatial convolutions \gian{like}
\beq\label{intop2}
B[\rho](x,t)=\ioma k(|y-x|)\,\rho(y,t)\,dy,
\eeq
with sufficiently regular kernels.

Optimal control problems of the above type often occur in industrial production processes. For instance,
consider a metallic workpiece consisting of two different component materials that tend to separate. Then a typical goal would be to monitor the production process in such a way that a desired distribution of the two materials (represented by the function $\rho_Q$) is realized during the time evolution in order to guarantee a wanted behavior of the workpiece; the deviation from the desired phase distribution is measured by the first summand in the 
cost $J$. The third summand of $J$ represents the costs due to the control action $u$; the size of the factors $\beta_i\ge 0$ then reflects the relative importance that the two conflicting interests ``realize the desired phase distribution as closely as possible'' and ``minimize the cost of the control action'' have for the manufacturer.

The state system (\ref{ss1})--(\ref{ss4}) is singular, with highly nonlinear
and nonstandard coupling. In particular, unpleasant nonlinear terms involving 
time derivatives occur in (\ref{ss1}),
and the expression $F'(\rho)$ in (\ref{ss2}) may become singular. Moreover,
the nonlocal term $B[\rho]$ is a source for possible analytical difficulties, and the
absence of the Laplacian in \eqref{ss2} may cause a low regularity of the 
order parameter $\rho$. We remark that the state system (\ref{ss1})--(\ref{ss4}) 
was recently analyzed in \cite{CGS3} for the case $u=0$ (no control); results 
concerning well-posedness and regularity were established.   

The mathematical literature on control problems for phase field systems involving equations
of viscous or nonviscous Cahn--Hilliard type is still scarce and quite recent. We refer in this connection to the works \cite{wn99,HW,CGS1,CGS2,CFGS1,CFGS2}. Control problems
for convective Cahn--Hilliard systems were studied in\pier{\cite{ZL1,ZL2}}, and a few
analytical contributions were made to the coupled Cahn--Hilliard/Navier--Stokes system
(cf. \cite{HW1,HW2,HW3,FRS}). The very recent contribution \cite{CGRS}
deals with the optimal control of a Cahn--Hilliard type system arising in the modeling of
solid tumor growth. 

The paper is organized as follows: in Section 2, we state the general assumptions and 
derive new regularity and stability results for the state system. In Section 3, we establish the
directional differentiability of the control-to-state operator, and the final Section 4
brings the main results of this paper, namely, the derivation of the first-order necessary
conditions of optimality.

Throughout this paper, we will use the following notation: we denote for a 
(real) Banach space $X$ by $\,\|\cdot\|_X$ its norm and the norm of 
$X\times X\times X$, by $X'$ its dual space, and by $\langle\cdot,\cdot\rangle_X$ 
the dual pairing between $X'$ and $X$. If $X$ is an inner product space, 
then the inner product is denoted by $(\cdot,\cdot)_X$. The only exception from this convention is given
by the $L^p$ spaces, $1\le p\le\infty$, for which we use the abbreviating notation
$\|\cdot\|_p$ \pier{for the norms.} Furthermore, 
we put 
$$H:=\lzo, \quad V:=H^1(\oma), \quad W:=\{w\in H^2(\oma) :\,\dn w=0\,\,\mbox{ a.\,e. on }\Gamma\}.$$
We have the dense and continuous embeddings 
$W\subset V\subset H\cong H'\subset V'\subset W'$, where $\langle u,v\rangle_V=(u,v)_H$ and
$\langle u,w\rangle_W=(u,w)_H$ for all $u\in H$, $v\in V$, and $w\in W$.
 
In the following, we will
make repeated use of Young's inequality
\begin{equation}
\label{young}
a\,b\le \delta \,a^2 + \frac 1{4\delta} b^2 \quad\mbox{for all }\,
a,b\in \rz\,\quad\mbox{and }\, \delta>0,
\end{equation}
as well as of the fact that for three dimensions of space and smooth domains 
the embeddings $\,V\subset L^p(\Omega)$, $1\le p\le 6$, and 
$\,H^2(\Omega)\subset C^0(\overline{\Omega})$ are continuous and 
(in the first case only for $1\le p<6$) 
compact. In particular, there are positive constants $\widetilde K_i$, $i=1,2,3$,
 which depend only on the domain $\oma$, such that
\beqa\label{embed1}
\|v\|_6&\!\!\le\!\!&\widetilde K_1\,\|v\|_V\quad\forall\,
v\in V,\\[1mm]
\label{embed2}
\|v\,w\|_H&\!\!\le\!\!&\|v\|_6\,\|w\|_3\,\le\,
\widetilde K_2\,\|v\|_V\,\|w\|_V\quad\forall\,v,w\in V, \qquad
\\[1mm]
\label{embed3}
\|v\|_\infty&\!\!\le\!\!&\widetilde K_3\,\|v\|_{H^2(\oma)}\quad\forall\,
v\in H^2(\oma).
\eeqa
\gianni{We also set for convenience
\beq
  Q_t := \Omega \times (0,t)
  \aand
  Q^t := \Omega \times (t,T) \pier{, \quad  \hbox{for }\, t\in (0,T).}
  \label{defQt}
\eeq
}\pier{Please note the difference between the subscript and the superscript in the above notation}.

\pier{About time derivatives of a time-dependent function $v$, we point out that we will use both the notations  $\pt v, \, \pt^2 v $ and the shorter ones $v_t, \, v_{tt} $.}


\section{Problem statement and results\\ for the state system}
\setcounter{equation}{0}
Consider the optimal control problem
(\ref{cost})--(\ref{uad}). 
We make the following assumptions on the data:

\vspace{5mm}
{\bf (A1)} \,\,$F=F_1+F_2$, where $F_1\in C^3(0,1)$ is convex, $F_2\in C^3[0,1]$,  and
\begin{equation}
\label{2.1}
\lim_{r\searrow 0} F_1'(r)=-\infty, \quad \lim_{r\nearrow 1} F_1'(r)= +\infty.
\end{equation}

{\bf (A2)} \,\,$\rho_0\in V$, $\pier{F'} (\rho_0) \in H$, $\mu_0\in W$, where
 $\,\mu_0\ge 0$\, a.\,e. in $\oma$,
\begin{equation}
\label{2.2}
\inf\,\{\rho_0(x): \,x\in\Omega\}>0, \quad \sup\,\{\rho_0(x):\, x\in\Omega\}<1\,. 
\end{equation}

{\bf (A3)} \,\,$g\in C^3[0,1]$ satisfies \pier{$g(\rho)\geq 0$ and} $g''(\rho)\le 0$ for all $\rho\in [0,1]$.

\vspace{3mm}
{\bf (A4)} \,\,The nonlocal operator \gianni{$B\colon L^1(Q)\to L^1(Q)$} satisfies the following conditions:

{\bf (i)} \,\,\,\,For every $t\in (0,T]$, we have
\begin{equation}\label{B1}
B[v]|_{Q_t}=B[w]|_{Q_t} \,\,\mbox{ whenever }\, v|_{Q_t}=w|_{Q_t}.
\end{equation}
{\bf (ii)} \,\,\,
\gianni{For all $p\in [2,+\infty]$, we have $B(L^p(Q_t))\subset L^p(Q_t)$ and
\beq
\label{B2}
 \|B[v]\|_{L^p(Q_t)}\,
\le C_{B,p}\left(1+\|v\|_{L^p(Q_t)}\right)
\eeq
for every $v\in L^p(Q)$ and $t\in(0,T]$.}

{\bf (iii)} \,\,For every $v,w\in L^1(0,T;H)$ and $t\in (0,T]$, it holds that
\beq
\label{Bsechs} 
\tint \|B[v](s)-B[w](s)\|_6\ds\,\le\,C_B\tint\|v(s)-w(s)\|_H\ds\,.
\eeq
{\bf (iv)} \,\,\,It holds, for every
$v\in L^2(0,T;V)$ and $t\in (0,T]$, that
\begin{align}
\label{B3}
&\|\nabla B[v]\|_{\lzht}\,\le\,C_B  \pier{\bigl( 1 + \,\|v\|_{\lzvt}\bigr)}.
\end{align}

\vspace{2mm}
{\bf (v)} \,\,\,\,For every $v\in H^1(0,T;\gianni H)$, we have $\,\pt B[v]\in L^2(Q)$
\gianni{and}
\beq
  \gianni{\|\pt B[v]\|_{L^2(Q)} \leq C_B \pier{\bigl( 1 + \|\pt v\|_{L^2(Q)}\bigr)}.  }
  \label{B4}
\eeq

{\bf (vi)} \,\,\,$B$ is continuously Fr\'echet differentiable as a mapping
from $L^2(Q)$ into $L^2(Q)$, and the Fr\'echet derivative $DB[\overline{v}]\in 
{\cal L}(L^2(Q),L^2(Q))$ of $B$ at $\overline{v}$ has for every $\overline{v}\in L^2(Q)$ 
and $t\in(0,T]$ the following properties: 
\begin{align}
\label{B5}
&\gianni{\|DB[\overline{v}](w)\|_{L^p(Q_t)}\,\le\,C_B\,\|w\|_{L^p(Q_t)} 
\quad\forall\, w\in L^p(Q),\quad \forall\, p\in[2,6],}\\[2mm]
\label{B6}
&\|\nabla(DB[\overline{v}](w))\|_{L^2(Q_t)} \,\le\,C_B\,\|w\|_{L^2(0,t;V)}\quad\forall \,
w\in L^2(0,T;V).
\end{align}

\vspace*{2mm}
In the above formulas, $C_{B,p}$ and $C_B$ denote given positive structural constants. 
\gianni{We also notice that \eqref{B5} implicitely requires that
$DB[\overline v](w)|_{Q_t}$ depends only on~$w|_{Q_t}$.
However, this is a consequence of~\eqref{B1}.}

\medskip

\gianni{%
The statements related to the control problem {\bf (CP)} depend on the assumptions made in the Introduction.
We recall them here.}

\smallskip
\gianni{{\bf (A5)}  \,\, $J$ and $\uad$ are defined by \eqref{cost} and \eqref{uad},
respectively, where
\begin{align}
\label{constCP}
&\beta_1\,,\,\beta_2\,,\,\beta_3\geq0, \quad \pier{\beta_1+ \beta_2 + \beta_3 > 0 ,  \aand R>0. }
\\[2mm]
\label{functCP}
&\rho_Q,\mu_Q\in\LQ2 , \quad
u_{\rm max} \in \LQ\infty
\aand
u_{\rm max} \geq0
\quad \hbox{a.e.\ in~$Q$.}
\end{align}
}%

{\sc Remark 1:} \,\,In view of \eqref{B5}, \gianni{for every $t\in [0,T]$ it holds that}
\beq\label{B8}
\|B[v]-B[w]\|_{L^2(Q_t)}\,\le\,C_B\,\|v-w\|_{L^2(Q_t)} \quad\forall\,v,w\in \lzq\,,
\eeq
that is, the condition (2.9) in \cite{CGS3} is fulfilled. Moreover, \eqref{B2} and 
\eqref{B3} imply that $B$ maps $L^2(0,T;V)$ into itself and that, for all \,$t\in (0,T]$\,
and \,$v\in L^2(0,T;V)$,
$$
\Bigl|\txinto \nabla B[v]\cdot\nabla v\dx\ds\Bigr|\,\le\,C_B\left(1\,+\,\|v\|_{\lzvt}^2
\right) ,
$$
which means that also the condition (2.10) in \cite{CGS3} is satisfied. Moreover, 
\pier{thanks to \eqref{B5} and \eqref{B6},}
there is some constant $\widetilde C_B>0$ such that
\beq\label{B9}
\|DB[\overline{v}](w)\|_{\lzvt}\,\le\,\widetilde{C}_B\,\|w\|_{\lzvt} \quad\forall\,
\overline{v}\in \lzq, \,\,\,\forall\,w \in L^2(0,T;V)\,.
\eeq

\vspace{3mm}
{\sc Remark 2:} \,\,We recall (cf.\ \cite{CGS3}) that
the integral operator \eqref{intop2} satisfies the conditions \eqref{B1} and \eqref{B2},
provided that
the integral
kernel $k$ belongs to $C^1(0,+\infty)$ and fulfills, with
suitable constants $C_1>0$, $C_2>0$, $0<\alpha<\frac 32$, $0<\beta<\frac 52$,  the growth conditions
$$
|k(r)|\le C_1\,r^{-\alpha}, \quad |k'(r)|\le C_2\,r^{-\beta}, \quad\forall\, r>0\,.
$$
In this case, we have $2\alpha<3$ and thus, for all $v,w\in L^1(0,T;H)$\, and $t\in (0,T]$,
\begin{align*}
&\tint\|B[v](s)-B[w](s)\|_6\ds\\[1mm]
&\le\,C_1\tint\Bigl(\ioma\Bigl|\ioma |y-x|^{-\alpha}
|v(y,s)-w(y,s)|\dy\Bigr|^6\dx\Bigr)^{1/6}ds\\[1mm]
&\le\,C_3\tint\Bigl(\ioma\Bigl|\Bigl(\ioma |y-x|^{-2\alpha}\dy\Bigr)^{1/2}\|v(s)-w(s)\|_H\Bigr|^6
\dx\Bigr)^{1/6}ds\\[1mm]
&\le\,C_4\tint\|v(s)-w(s)\|_H\ds,
\end{align*}
with global constants $C_i$, $3\le i\le 4$; the condition \eqref{Bsechs} is thus satisfied. Also
condition \eqref{B3} holds true in this case: indeed, for every $t\in (0,T]$ and $v\in L^2(0,T;V)$,
we find, since $\frac {6\beta}5<3$, that
\begin{align*}
&\|\nabla B[v]\|^2_{\lzht}\,\le\,C_2\txinto\Bigl|\ioma|y-x|^{-\beta}|v(y,s)|\dy\Bigr|^2\dx\ds\\[1mm]
&\le\,C_5\txinto\Bigl(\ioma|y-x|^{-6\beta/5}dy\Bigr)^{5/3}\|v(s)\|_6^2\dx\ds\\[1mm]
&\le \,C_6\tint\|v(s)\|_V^2\ds\,.
\end{align*}   
Finally, since the operator $B$ is linear in this case, we have 
$DB[\overline{v}]=B$ for every $\overline{v}\in\lzq$, and thus also {\bf (A4)(v)} and
\eqref{B5}--\eqref{B9} are fulfilled. 
 Notice that the above growth conditions are 
met by, e.\,g., the three-dimensional Newtonian potential, where $k(r)=\frac c r$ with some
$c\not =0$.  

\vspace{2mm} We also note that {\bf (A2)} implies $\mu_0\in C(\overline\Omega)$, and \pier{{\bf (A1)} and \eqref{2.2} ensure that both $F(\rho_0)$ and  $F'(\rho_0)$ are in $L^\infty (\Omega)$, whence in $H$.}
Moreover, the logarithmic potential \eqref{logpot}
obviously fulfills the condition \eqref{2.1} in {\bf (A1)}.  

We have the following existence and regularity result for the state system.

\vspace{5mm}
{\sc Theorem 2.1:} \quad {\em Suppose that} {\bf (A1)}--{\bf (A\gianni5)} {\em are 
satisfied. Then the state system} (\ref{ss1})--(\ref{ss4}) {\em has for every $u\in\uad$ a unique solution}
$(\rho,\mu)$ {\em such that}
\begin{eqnarray}
\label{reguss1}
&&\rho\in H^2(0,T;H)\cap \pier{W^{1,\infty}(0,T;\lio) \cap  H^1(0,T;V)} ,\\[1mm]
\label{reguss3}
&&\mu\in W^{1,\infty}(0,T;H)\cap H^1(0,T;V)\cap L^\infty(0,T;W) \mathrel{\gianni\subset} L^\infty(Q).\qquad
\end{eqnarray}   
{\em Moreover, there are constants $0<\rho_*<\rho^*<1$, $\mu^*>0$, and $K_1^*>0$, 
which depend only on the given data, such that for every $u\in\uad$ the 
corresponding solution $(\rho,\mu)$ satisfies}
\begin{eqnarray}
\label{ssbounds1}
&& 0<\rho_*\le \rho\le\rho^*<1,\,\quad 0\le\mu\le\mu^*, 
\,\quad\mbox{a.\,e. in }\,Q,\\[1mm]
 \label{ssbounds2}
&&\|\mu\|_
{W^{1,\infty}(0,T;H)\cap H^1(0,T;V)\cap L^\infty(0,T;W)\cap
  L^\infty(Q)}\nonumber\\[1mm]
&&+\,\|\rho\|_{H^2(0,T;H)\cap W^{1,\infty}(0,T;\lio)\cap H^1(0,T;V)}\,\le\,K^*_1.
\end{eqnarray}

\vspace{2mm}
{\sc Proof:} \,\,\,In the following, we denote by $C_i>0$, $i\in\nz$, constants which
may depend on the data of the control problem {\bf (CP)} but not on the special
choice of $u\in\uad$. First, we note that
in \cite[Thms.\,2.1,\,\,2.2]{CGS3} it has been shown that under the given assumptions
there exists for $u\equiv 0$ a unique solution $(\rho,\mu)$ with the properties
\begin{align}\label{zw1}
&0<\rho<1,\quad \mu\ge 0, \quad\mbox{a.\,e. in }\,Q,\\[2mm]
\label{zw2}
&\rho\in L^\infty(0,T;V), \quad \pt\rho\in L^6(Q), \\[2mm]
\label{zw3}
& \mu\in H^1(0,T;H)\cap L^\infty(0,T;V)\cap L^\infty(Q)\cap L^2(0,T;W^{2,3/2}(\oma)).  
\end{align}

A closer inspection of the proofs in \cite{CGS3} reveals that the line of argumentation used there
(in particular, the proof that $\mu$ is nonnegative) carries over with only minor modifications to general right-hand sides $u\in\uad$. We thus infer
that \eqref{ss1}--\eqref{ss4} enjoys for every $u\in\uad$ a unique solution satisfying \eqref{zw1}--\eqref{zw3};
more precisely, there is some $C_1>0$ such that
\begin{align}
\label{bd1}
&\|\mu\|_{H^1(0,T;H)\cap L^\infty(0,T;V)\cap L^\infty(Q)\cap L^2(0,T;W^{2,3/2}(\oma))}\nonumber\\[1mm]
&+ \,\|\rho\|_{L^\infty(0,T;V)}\,+\,\|\pt\rho\|_{L^6(Q)}\,\le\,C_1 \quad\forall\,u\in\uad.
\end{align}
Moreover, invoking \eqref{zw1}, and \eqref{B2} for $p=+\infty$, we find that
$$
\|B[\rho]\|_{L^\infty(Q)}\,\le\,C_2 \quad\forall\,u\in\uad,
$$
and it follows from \eqref{bd1} and the general assumptions on $\rho_0$, $g$, and $F$, that there are constants
$\rho_*,\rho^*$ such that, for every $u\in\uad$,
\begin{eqnarray*}
&&0<\rho_*\,\le\,\inf\,\{\rho_0(x):x\in\oma\}\,\le\, \sup\,\{\rho_0(x):x\in\oma\}\,\le\,\rho^*<1,\\[1mm]
&&F'(\rho)+B[\rho]-\mu\,g'(\rho)\,\le\,0 \,\,\mbox{ if \,$0 <\rho\le \rho_*$},\\[1mm]
&&F'(\rho)+B[\rho]-\mu\,g'(\rho)\,\ge\,0 \,\,\mbox{ if \,$\rho^*\le\rho<1$.}
\end{eqnarray*}
Therefore, multiplying \eqref{ss2} by  the positive part $(\rho-\rho^*)^+$ of $\rho-\rho^*$, and integrating
over~$Q$, we find that
\begin{align*}
0&=\int_0^T\!\!\!\ioma \pt\rho\,(\rho-\rho^*)^+\dx\dt\, +\int_0^T\!\!\!\ioma (F'(\rho)+B[\rho]-\mu\,g'(\rho))
(\rho-\rho^*)^+\dx\dt\\[1mm]
&\ge \frac 12\ioma \left|(\rho(t)-\rho^*)^+\right|^2\dx,
\end{align*} 
whence we conclude that $\,(\rho-\rho^*)^+=0$, and thus $\rho\le\rho^*$, almost everywhere in~$Q$. The other
bound for $\rho$ in \eqref{ssbounds1} is proved similarly.

It remains to show the missing bounds in \eqref{ssbounds2} (which then also imply
the missing regularity claimed in \eqref{reguss1}--\eqref{reguss3}). To this end, we
employ a bootstrapping argument.

First, notice that {\bf (A3)} and the already proved 
bounds \eqref{bd1} and \eqref{ssbounds1} imply that the expressions
$\,\mu\,g'(\rho)\,\pt\rho$ \,and\,$(1+2g(\rho))\,\pt\mu\,$ are bounded in~$\lzq$. 
Hence, by comparison in
\eqref{ss1},
the same holds true for~$\Delta\mu\,$, and thus standard elliptic estimates yield that
\beq
\label{bd2}
\|\mu\|_{L^2(0,T;\gianni W)}\,\le\,C_3 \quad\forall\,u\in\uad.
\eeq
Next, observe that {\bf (A1)} and \eqref{ssbounds1} imply that $\,\|F'(\rho)\|_{L^\infty(Q)}\,\le\,C_4$,
and comparison in \eqref{ss2}, using {\bf (A3)},  yields that
\beq
\label{bd3}
\|\pt\rho\|_{L^\infty(Q)}\,\le\,C_5 \quad\forall\,u\in\uad.
\eeq
In addition, we infer from the estimates shown above, and using \eqref{B3}, that the
right-hand side of the identity
\begin{align}
\label{nablarhot}
\nabla\rho_t = -F''(\rho)\,\nabla\rho-\nabla B[\rho]+g'(\rho)\,\nabla\mu
+\mu\, g''(\rho)\,\nabla\rho
\end{align}
is bounded in $\lzq$, so that
\beq \label{bd4}
\|\pt\rho\|_{L^2(0,T;V)}\,\le\,C_6 \quad\forall\,u\in\uad.
\eeq
We also note that the time derivative $\,\pt(-F'(\rho)-B[\rho]+\mu g'(\rho))\,$ exists and is bounded
in~$\lzq$ \gianni{(cf.~\eqref{B4})}. We thus have
\beq \label{bd5}
\|\rho_{tt}\|_{\lzq}\,\le\,C_7 \quad\forall\,u\in\uad.
\eeq
At this point, we observe that Eq.~\eqref{ss1} can be written in the form
$$
a\,\pt\mu+\mu\,\pt a-\Delta\mu=b,\quad\mbox{with }\,a:=1+2g(\rho), \quad b:=u+\mu\,g'(\rho)\,\pt\rho,
$$
where,  thanks to the above estimates, we have, for every $u\in\uad$, 
\begin{align}
\label{bd5bis}
&\pier{\left\| a\right\|_{L^\infty (Q)} +  \left\|\pt a\right\|_{L^\infty (Q)}  +  \left\| b \right\|_{L^\infty (Q)} \leq  C_8}, \\[1mm]
\label{bd6}
&\bigl\|\pt^2 a\bigr\|_{\lzq}\,=\,2\,\|g''(\rho)\rho_t^2+g'(\rho)\rho_{tt}\|_{\lzq}
\,\le\,\pier{C_9},\\[1mm]
\label{bd7}
&\|\pt b\|_{\lzq}\,=\,\|u_t+ \mu_t g'(\rho)\rho_t+\mu g''(\rho)\rho_t^2+\mu g'(\rho)\rho_{tt}\|
_{\lzq}\,\le\,\pier{C_{10}}.
\end{align}
Since also $\mu_0\in W$, we are thus in the situation of \cite[Thm.\,3.4]{CGSRendiconti}, whence
we obtain that $\pt\mu \in L^\infty(0,T;H)\cap L^2(0,T;V)$ and $\mu\in L^\infty(0,T;W)$. Moreover,
a closer look at the proof of \cite[Thm.\,3.4]{CGSRendiconti} reveals that we also have the estimates
\begin{align}
\label{muW}
&\|\pt\mu\|_{L^\infty(0,T;H)\cap L^2(0,T;V)}\,+\,\|\mu\|_{L^\infty(0,T;W)}\,\le\,\pier{C_{11}}.
\end{align}
This concludes the proof of the assertion.
\qed

\vspace{5mm}
{\sc Remark 3:}
\,\,\,From the estimates \eqref{ssbounds1} and \eqref{ssbounds2}, and using 
the continuity of the
embedding $V\subset L^6(\Omega)$, we can without loss of generality (by possibly choosing a larger $K_1^*$) assume 
that also
\begin{align}
\label{ssbounds3}
&\nonumber\max_{0\le i\le 3}\,\|F^{(i)}(\rho)\|_{L^\infty(Q)}\,+\,\max_{0\le i\le 3}
\|g^{(i)}(\rho)\|_{L^\infty(Q)}\nonumber\\[1mm]
&+\,\|\nabla\mu\|_{L^\infty(0,T;L^6(\Omega)^3)}\,+\,\|\pt\mu\|_{L^2(0,T;V)}\nonumber\\[1mm]
&+\,\|B[\rho]\|_{H^1(0,T;\lzo)\cap L^\infty(Q)\cap L^2(0,T;V)} \le \,K_1^* \quad\forall\,u\in\uad\,.
\end{align}

\vspace{5mm}
According to Theorem 2.1, the control-to-state mapping $\,{\cal S}:\uad\ni u\mapsto (\rho,\mu)$\, is well defined. We now study its stability properties. We have the following result.

\vspace{5mm}
{\sc Theorem 2.2:} \quad {\em Suppose that} {\bf (A1)}--{\bf (A\gianni5)} {\em are fulfilled,
and let $u_i\in\uad$, $i=1,2$, be given and $(\rho_i,\mu_i)={\cal S}(u_i)$, $i=1,2$, be
the associated solutions to the state system}  \eqref{ss1}--\eqref{ss4}. 
{\em Then there
exists a contant $K_2^*>0$, which depends only on the data of the problem, such that, for
every $t\in (0,T]$,}
\beqa
\label{stabu}
&&\|\rho_1-\rho_2\|_{H^1(0,t;H)\cap L^\infty(0,t;L^6(\oma))}
\,+\,\|\mu_1-\mu_2\|_{H^1(0,t;H)\cap L^\infty(0,t;V)\cap L^2(0,t;\gianni W)}\nonumber\\[2mm]
&&\le\,K_2^*\,\|u_1-u_2\|_{\lzht}\,.
\eeqa

\vspace{2mm}
{\sc Proof:} \,\,Taking the difference of the equations satisfied by $(\rho_i,\mu_i)$,
$i=1,2$, and setting $u:=u_1-u_2$, $\rho:=\rho_1-\rho_2$,
$\mu:=\mu_1-\mu_2$, we first observe that we have almost everywhere in $Q$ the identities
\begin{align}\label{diff1}
&(1+2g(\rho_1))\,\pt\mu \,+\,g'(\rho_1)\,\pt\rho_1\,\mu -\Delta\mu
\,+\,2(g(\rho_1)-g(\rho_2))\,\pt\mu_{2}\nonumber\\[1mm]
&\quad=u\,-\,(g'(\rho_1)-g'(\rho_2))\,\pt\rho_1\,\mu_2\,-\,g'(\rho_2)\,\mu_2\,\pt\rho\,,\\[3mm]
\label{diff2}
&\pt\rho\,+\,F'(\rho_1)\,-F'(\rho_2)\,+\,B[\rho_1]-B[\rho_2]\nonumber\\[1mm]
&\quad=g'(\rho_1)\,\mu\,+\,(g'(\rho_1)-g'(\rho_2))\,\mu_2\,,
\end{align}
\pier{as well as}
\begin{align}
\label{diff4}
\dn\mu=0\quad\mbox{a.\,e. on }\,\Sigma, \quad \mu(\cdot,0)=\rho(\cdot,0)=0 \quad
\mbox{a.\,e. in }\,\oma.
\end{align}

Let $t\in (0,T]$ be arbitrary. In the following, we repeatedly use the global estimates
\eqref{ssbounds1}, \gianni{\eqref{ssbounds2}}, 
and \eqref{ssbounds3}, without further reference. Moreover, we denote by
$c>0$ constants that may depend on the given data of the state system, but not on 
the choice of $u_1,u_2\in\uad$; the meaning of $c$ 
may change between and even within lines. We establish the validity of \eqref{stabu}
in a series of steps.

\vspace{2mm}
\underline{\sc Step 1:} \,\,\,
To begin with, we first 
observe that $$\,(1+2g(\rho_1))\mu\pt\mu+g'(\rho_1)\pt\rho_1\mu^2=\pt\Bigl(
\Bigl(\frac 12+g(\rho_1)\Bigr)\mu^2\Bigr).$$ Hence, multiplying \eqref{diff1} by $\mu$ and
integrating over $Q_t$ and by parts, we obtain that
\beq
\label{p221}
\ioma \Big(\frac 12+g(\rho_1(t))\Big)\mu^2(t)\dx\,+\txinto|\nabla\mu|^2\dx\ds\, \le\,\sum_{j=1}^3
|I_j|\,,
\eeq
where the expressions $I_j$, $j=1,2,3$, defined below, are estimated as follows: first, we apply {\bf (A3)}, the
mean value theorem, and H\"older's and Young's \pier{inequalities}, to find, for every $\gamma>0$ (to be chosen later), that
\begin{align}
\label{p222}
&I_1:=-2\!\txinto\!\!(g(\rho_1)-g(\rho_2))\,\pt\mu_2\,\mu\dx\ds\le c
\!\! \tint\!\!\|\pt\mu_2(s	)\|_6
\,\|\mu(s)\|_3\,\|\rho(s)\|_2\ds\nonumber\\[1mm]
&\hspace*{4mm} \le\,\gamma\tint\|\mu(s)\|_V^2\ds\,+\,\frac c\gamma\tint
\|\pt\mu_2(s)\|_V^2\,\|\rho(s)\|_H^2 \ds,
\end{align}  
where it follows from \eqref{ssbounds3} that the mapping $\,\,s
\mapsto\|\pt\mu_2(s)\|_V^2\,\,$ belongs to $L^1(0,T)$.
Next, we see that 
\begin{align}
\label{p223}
&I_2:=\txinto \Bigl(u-(g'(\rho_1)-g'(\rho_2))\pt\rho_1\,\mu_2\Bigr)\,\mu\dx\ds\nonumber\\[1mm]
&\hspace*{5mm}\le c\txinto(|u|+|\rho|)|\mu|\dx\ds
\,\le\,c\txinto (u^2+\rho^2+\mu^2)\dx\ds\,.
\end{align}
Finally, Young's inequality yields that 
\beq
\label{p224}
I_3:=-\txinto g'(\rho_2)\,\mu_2\,\rho_t\,\mu\dx\ds\,\le\,\gamma\txinto\rho_t^2\dx\ds\,+\,
\frac c \gamma \txinto\mu^2\dx\ds\,.
\eeq
Combining \eqref{p221}--\eqref{p224}, and recalling that $g(\rho_1)$ is nonnegative,
we have thus shown the estimate
\begin{align}
\label{p225}
&\frac 1 2\,\|\mu(t)\|_H^2\,+\,(1-\gamma)\tint\|\mu(s)\|_V^2\ds\,\le\,\gamma\txinto\rho_t^2\dx\ds\,+\,c \txinto u^2\dx\ds\nonumber\\[1mm]
&+\, c\left(1+\gamma^{-1}\right)\tint\bigl(\|\mu(s)\|_H^2\,+\,(1+\|\pt\mu_2(s)\|_V^2)
\,\|\rho(s)\|_H^2\bigr)\ds\,.
\end{align}

Next, we add $\rho$ on both sides of \eqref{diff2} and multiply the result by $\rho_t$.
Integrating over $Q_t$, using \pier{the Lipschitz continuity of $F'$ (when restricted to 
$[\rho_*, \rho^*]$),  \eqref{B8} and Young's inequality,} we easily find the estimate
\beq
\label{p226}
(1-\gamma)\txinto \rho_t^2\dx\ds\,+\,\frac 12 \,\|\rho(t)\|_H^2\,\le\,\frac c \gamma
\txinto(\rho^2+\mu^2)\dx\ds\,.
\eeq
Therefore, combining \eqref{p225} with \eqref{p226}, choosing $\gamma>0$ small enough,
and invoking Gronwall's lemma, we have shown that
\beq
\label{p227}
\|\mu\|_{L^\infty(0,t;\gianni H)\cap \lzvt}\,+\,\|\rho\|_{H^1(0,t;\gianni H)}\,\le\,c
\|u\|_{\lzht} \quad\forall \,t\in(0,T].
\eeq

\vspace{3mm}
\underline{\sc Step 2:}  \,\,\,Next, we multiply \eqref{diff2} by $\,\rho|\rho|\,$  and integrate
over $Q_t$. We obtain 
\beq
\label{p228}
\frac  13 \,\|\rho(t)\|_3^3 \,\le\, \sum_{j=1}^3 |J_j|, 
\eeq
where the expressions $J_j$, $1\le j\le 3$, are estimated as follows: at first, we \pier{simply have}
\begin{align}
\label{p229}
J_1:&= \txinto (-F'(\rho_1)+F'(\rho_2)+\mu_2(g'(\rho_1)-g'(\rho_2)))\,\rho\,|\rho|\dx\ds\nonumber\\[1mm]
&\le c\tint\|\rho(s)\|_3^3\ds.
\end{align}
Moreover, invoking \eqref{p227}, H\"older's inequality, as well as  the global bounds once more,
\begin{align}
\label{p2210}
J_2:&= \txinto\mu\,g'(\rho_1)\,\rho\,|\rho|\dx\ds\,\le\,c\tint\|\mu(s)\|_6\,\|\rho(s)\|_2\,\|\rho(s)\|_3\ds
\nonumber\\[1mm]
&\le\tint \|\rho(s)\|_3^3\ds\,+\,c\tint\|\mu(s)\|_V^{3/2}\,\|\rho(s)\|_H^{3/2}\ds\nonumber\\[1mm]
&\gianni{\le\tint\|\rho(s)\|_3^3\ds\,+\,c\,\|\rho\|_{L^\infty(0,t;H)}^{3/2}\,\|\mu\|_{L^{3/2}(0,t;V)}^{3/2}}\nonumber\\[1mm]
&\le\tint\|\rho(s)\|_3^3\ds\,+\,c\,\|\rho\|_{L^\infty(0,t;H)}^{3/2}\,\|\mu\|_{\lzvt}^{3/2}\nonumber\\[1mm]
&\le \tint\|\rho(s)\|^3\ds\,+\,c\,\|u\|_{\lzht}^3\,.
\end{align}
In addition, condition \eqref{Bsechs}, H\"older's inequality, and \eqref{p227}, yield that
\begin{align}
\label{p2211}
J_3:&= -\txinto(B[\rho_1]- B[\rho_2])\,\rho\,|\rho|\dx\ds\,\nonumber\\[1mm]
&\le\,c\tint\pier{\|\rho(s)\|_3\,\|\rho(s)\|_2 \, \|B[\rho_1](s)-B[\rho_2](s)\|_6 }\ds \nonumber\\[1mm]
&\le\,c\,\sup_{0\le s\le t}\|\rho(s)\|_3\,\|\rho\|_{L^\infty(0,t;H)}
\tint\|\rho(s)\|_H\ds\nonumber\\[1mm]
&\le\,\frac 16\,\sup_{0\le s\le t}\|\rho(s)\|_3^3\,+\,c\,\|u\|^3_{\lzht}\,.
\end{align}
Combining the estimates \eqref{p228}--\eqref{p2211}, and invoking Gronwall's lemma, we can easily infer
that
\beq
\label{p2212}
\|\rho\|_{L^\infty(0,t;L^3(\oma))}\,\le\,c\,\|u\|_{\lzht} \quad\mbox{for all }\,t\in (0,T]. 
\eeq

\vspace{3mm}
\underline{Step 3:} \,\,\,With the above estimates proved, the road is paved for multiplying \eqref{diff1}
by $\mu_t$. Integration over $Q_t$ yields that
\begin{align}
\label{p2213}
\txinto (1+2g(\rho_1))\,\mu_t^2\dx\ds\,+\,\frac 12\,\|\nabla\mu(t)\|_H^2
\,\le\,\sum_{j=1}^5\,|K_j|,
\end{align}
where the expressions $K_j$, $1\le j\le 5$, are estimated as follows: at first, using the global bounds and
Young's inequality, we have for every $\gamma>0$ (to be specified later) the bound
\begin{align}
\label{p2214}
K_1:&=\,-\txinto g'(\rho_1)\,\pt\rho_1\,\mu\,\mu_t\dx\ds \,\le\,\gamma\txinto\mu_t^2\dx\ds
\,+\,\frac c\gamma\txinto\mu^2\dx\ds\nonumber\\[1mm]                           
&\le \,\gamma\txinto \mu_t^2\dx\ds\,+\,\frac c\gamma\,\|u\|^2_{\lzht}\,.
\end{align}
Next, thanks to the mean value theorem, and employing \eqref{ssbounds3} and \eqref{p2212}, we find that
\begin{align}
\label{p2215}
K_2:&=\,-\,2\txinto(g(\rho_1)-g(\rho_2))\,\pt\mu_2\,\mu_t\dx\ds\,\le\,c\txinto |\rho|\,|\pt\mu_2|\,|\mu_t|\dx\ds
\nonumber\\[1mm]
&\le\,c\tint\|\rho(s)\|_3\,\|\pt\mu_2(s)\|_6\,\|\mu_t(s)\|_2\ds\nonumber\\[1mm]
&\le\,\gamma\txinto \mu_t^2\dx\ds\,+\,\frac c\gamma\,\|\rho\|_{L^\infty(0,t;L^3(\Omega))}^2\tint\|\pt\mu_2(s)\|_V^2\ds
\nonumber\\[1mm]
&\le\,\gamma\txinto\mu_t^2\ds\,+\,\frac c\gamma\,\|u\|^2_{\lzht}\,.
\end{align}
Moreover, \pier{we infer that}
\begin{align}
\label{p2216}
K_3:=\txinto u\,\mu_t\dx\ds\,\le\,\gamma\txinto\mu_t^2\dx\ds\,+\,\frac c\gamma\,\|u\|^2_{\lzht}\,,
\end{align}
as well as, invoking the mean value theorem once more,  
\begin{align}
\label{p2217}
K_4:&=\,-\txinto (g'(\rho_1)-g'(\rho_2))\,\pt\rho_1\,\mu_2\,\mu_t\dx\ds\,\le\,c\txinto| \rho|\,|\mu_t|\dx\ds
\nonumber\\[1mm]
&\le\,\gamma\txinto \mu_t^2\dx\ds\,+\,\frac c\gamma\,\|u\|^2_{\lzht}\,,
\end{align}
and, finally, using \eqref{p227} and Young's inequality,
\begin{align}
\label{p2218}
K_5:&=\,-\txinto g'(\rho_2)\,\mu_2\,\rho_t\,\mu_t\dx\ds\,\le\,c\txinto|\rho_t|\,|\mu_t|\dx\ds\nonumber\\[1mm]
&\le \,\gamma\txinto\mu_t^2\dx\ds\,+\,\frac c\gamma\,\|\rho\|^2_{H^1(0,t;H)}\nonumber\\[1mm]
&\le\,\gamma\txinto\mu_t^2\dx\ds\,+\,\frac c\gamma\,\|u\|^2_{\lzht}\,.
\end{align}

Now we combine the estimates \eqref{p2213}--\eqref{p2218} and choose $\gamma>0$ appropriately
small. It then follows that
\begin{align}
\label{p2219}
\|\mu\|_{H^1(0,t;H)\cap L^\infty(0,t;V)}\,\le \,c\,\|u\|_{\lzht}\,.
\end{align}

Finally, we \gianni{come back to \eqref{diff1} and}
employ the global bounds \eqref{ssbounds1}, \eqref{ssbounds2}, \eqref{ssbounds3}, and the estimates shown above, to conclude that
\begin{align}
\|\Delta\mu\|_{\lzht}\,&\le\,c\left(\|\mu_t\|_{\lzht}\,+\,\|\mu\|_{\lzht}\,+\,\|\rho_t\|_{\lzht}
\right.\nonumber\\[1mm]
&\left.\hspace*{10mm} \pier{{}+\,\|\rho\|_{\lzht}}\,+\,\|u\|_{\lzht}\right)\,+\,c\,\|\rho\,\pt\mu_2\|_{\lzht}
\nonumber\\[1mm]
&\le\,c\,\|u\|_{\lzht}\,,
\end{align}
where the last summand on the right-hand side was estimated as follows:
\begin{align*}
&\txinto\,|\rho|^2\,|\pt\mu_2|^2\dx\ds\,\le\,c\tint\|\pt\mu_2(s)\|_6^2\,\|\rho(s)\|_3^2\ds\\[1mm]
&\le\,c\,\|\rho\|^2_{L^\infty(0,t;L^3(\oma))}\tint\|\pt\mu_2(s)\|_V^2\ds\,\le\,c\,\|u\|^2_{\lzht}\,.
\end{align*}

Invoking standard elliptic estimates, we have thus shown that
\beq
\label{p2220}
\|\mu\|_{L^2(0,t;\gianni W)}\,\le\,c\,\|u\|_{\lzht}\,.
\eeq

\vspace{3mm}
\underline{\sc Step 4:} \,\,\,It remains to show the $L^\infty(0,t;L^6(\oma))$--\,estimate for $\rho$.
To this end, we multiply \eqref{diff2} by $\,\rho|\rho|^4\,$ and integrate over $Q_t$. It follows
that
\beq
\label{p2221}
\frac 16\,\|\rho(t)\|_6^6\,\le\,\sum _{j=1}^3 |L_j|,
\eeq
where quantities $L_j$, $1\le j\le 3$, are estimated as follows: at first, we obtain from the
global estimates \eqref{ssbounds2} and \eqref{ssbounds3}, that
\begin{align}
\label{p2222}
L_1:&=\txinto\bigl(-F'(\rho_1)+F'(\rho_2)+\mu_2(g'(\rho_1)-g'(\rho_2))\bigr)\,\rho\,|\rho|^4\dx\ds
\nonumber\\[1mm]
&\le\,c\tint\|\rho(s)\|_6^6\ds\,.
\end{align}
Moreover, from \eqref{p2219} and H\"older's and Young's inequalities \pier{we conclude that}
\begin{align}
\label{p2223}
L_2:&=\txinto g'(\rho_1)\,\mu\,\rho\,|\rho|^4\dx\ds\,\le\,c\tint\|\mu(s)\|_6\,\|\rho(s)\|_6^5\ds
\nonumber\\[1mm]
&\gianni{\le\,c\,\|\mu\|_{L^\infty(0,t;V)}\,\, \|\rho\|_{L^5(0,t;L^6(\Omega))}^5
\leq \,c\,\|\mu\|_{L^\infty(0,t;V)}^6 + c\, \|\rho\|_{L^5(0,t;L^6(\Omega))}^6}
\nonumber\\[1mm]
&\,\le\,c\,\|u\|^6_{\lzht}\,+\,c\tint\|\rho(s)\|_6^6\ds\,.
\end{align}
Finally, we employ \eqref{Bsechs} and \eqref{p227} to infer that
\begin{align}
\label{p2224}
L_3:&=\,-\txinto(B[\rho_1]-B[\rho_2])\,\rho\,|\rho|^4\dx\ds\nonumber\\[1mm]
&\le\,c\tint\|B(\rho_1](s)-B[\rho_2](s)\|_6\,\|\rho(s)\|_6^5\ds\nonumber\\[1mm]
&\le\,c\,\sup_{0\le s\le t}\|\rho(s)\|_6^5\tint\|\rho(s)\|_H\ds\nonumber\\[1mm]
&\le\,\frac 1{12}\,\sup_{0\le s\le t}\|\rho(s)\|_6^6 \,+\,c\,\|u\|^6_{\lzht}\,.
\end{align} 
Combining the estimates \eqref{p2221}--\eqref{p2224}, and invoking Gronwall's lemma,
\pier{then we} readily find the estimate
$$
\|\rho\|_{L^\infty(0,t;L^6(\oma))}\,\le\,c\,\|u\|_{\lzht}\,,
$$
which concludes the proof of the assertion.\qed

\section{Directional differentiability\\ of the control-to-state mapping}
\setcounter{equation}{0}

In this section, we prove the relevant differentiability properties of the
solution operator $\cs$. To this end, 
we introduce the spaces 
\begin{align*}
&\cx:=H^1(0,T;H)\cap \liq,\\
& \cy:=H^1(0,T;H)\times \left(L^\infty (0,T;H)
\cap L^2(0,T;V)\right),
\end{align*}
endowed with their natural norms
\begin{align*}
\|u\|_{\cx}&:=\,\|u\|_{H^1(0,T;H)}\,+\,\|u\|_{\liq} \quad\forall\,u\in\cx,\\[1mm]
\|(\rho,\mu)\|_{\cy}&:=\|\rho\|_{H^1(0,T;H)}\,+\,\|\mu\|_{L^\infty(0,T;H)}
\,+\,\|\mu\|_{L^2(0,T;V)} \quad\forall\,(\rho,\mu)\in\cy\,,
\end{align*}
and consider the control-to-state operator $\cs$ as a mapping between $\uad\subset\cx$
and $\cy$. Now let $\bu\in\uad$ be fixed and
put $(\br,\bmu):=\cs(\bu)$. 
We then study the linearization of the state system
\eqref{ss1}--\eqref{ss4} at~$\bu$, \gianni{which is given~by}:
\begin{align}
\label{ls1}
&(1+2g(\br))\,\eta_t\,+\,2g'(\br)\,\bmu_t\,\xi\,+\,g'(\br)\,\br_t\,\eta 
\,+\,\bmu\, g''(\br)\,\br_t\,\xi
\,+\,\bmu g'(\br)\,\xi_t\nonumber\\ 
&\quad\quad-\Delta\eta = h \quad\mbox{a.\,e.\ in }\,Q,\\[1mm]
\label{ls2}
&\xi_t+F''(\br)\,\xi + DB[\br](\xi)=\bmu\, g''(\br)\,\xi\,+\,g'(\br)\,\eta
\quad\mbox{a.\,e.\ in }\,Q,\\[1mm]
\label{ls3}
&\hspace*{3cm} \dn\eta=0 \quad\mbox{a.\,e.\ on }\, \Sigma,\\[1mm]
\label{ls4}
&\hspace*{2.9cm} \eta(0)=\xi(0)=0 \quad\mbox{a.\,e.\ in }\,\oma.
\end{align}

Here, $h\in\cx$ \pier{must satisfy} $\bu+\gianni\blam h\in\uad$ for some $\gianni\blam>0$.
Provided that the system \eqref{ls1}--\eqref{ls4} has for any such $h$ 
a unique solution pair $(\xi,\eta)$, we expect that the directional
derivative $\delta\cs(\bu;h)$ of $\cs$ at $\bu$ in the direction $h$ (if it exists) ought to 
be given by $(\xi,\eta)$. 
\gianni{In fact, the above problem makes \juerg{sense} for every $h\in\LQ2$, 
and it is uniquely solvable under this weaker assumption.}

\vspace{5mm}
{\sc Theorem 3.1:} \quad {\em Suppose that the general hypotheses} {\bf (A1)}--{\bf (A5)}
{\em are satisfied and let $h\in\LQ2$.
Then, the linearized problem} \accorpa{ls1}{ls4} 
{\em has a unique solution $(\xi,\eta)$ satisfying
\Bsist
  && \xi \in \H1H \cap \L\infty\Lsei,
  \label{regxi}
  \\
  && \eta \in \H1H \cap \L\infty V \cap \L2W.
  \label{regeta}
\Esist
}%

\vspace{3mm}
{\sc Proof:} 
We first prove uniqueness.
Since the problem is linear,
we take $h=0$ and show that $(\xi,\eta)=(0,0)$.
We add $\eta$ and $\xi$ to both sides of equations \eqref{ls1} and~\eqref{ls2}, respectively,
then multiply by $\eta$ and~$\xi_t$, integrate over~$Q_t$, and sum~up.
By observing that
\Beq
  (1+2g(\br)) \eta_t \eta + g'(\br) \br_t |\eta|^2
  = \pt \bigl[ \juerg{\bigl({\smfrac 12} + g(\br)\bigr)} |\eta|^2 \bigr],
  \non
\Eeq
and recalling that $g\geq0$, we obtain
\Beq
  \frac 12 \iO |\eta(t)|^2 \dx
  + \iot \normaV{\eta(s)}^2 \ds
  + \frac 12 \iO |\xi(t)|^2 \dx
  + \intQt |\xi_t|^2 \dx\ds
  \leq \sum_{j=1}^3 H_j,
  \non
\Eeq
where the terms $H_j$ are defined and estimated as follows. 
We have
\Bsist
  && H_1 := - \intQt 2\, \pier{g'(\br)}\,\bmu_t\, \xi \,\eta\dx\ds
  \leq c \iot \norma{\bmu_t(s)}_3 \, \norma{\xi(s)}_2 \, \norma{\eta(s)}_6 \ds
  \non
  \\
  && \leq \frac 12 \iot \normaV{\eta(s)}^2 \ds
  + c \iot \normaV{\bmu_t(s)}^2 \, \juerg{\norma{\xi(s)}^2_2} \ds\,,
  \non
\Esist
and we notice that the function $s\mapsto\normaV{\bmu_t(s)}^2$ belongs to $L^1(0,T)$,
 by~\eqref{muW} for~$\bmu$.
Next, we easily have the estimate
\Bsist
  && H_2 := \intQt \bigl(
    \eta
    - \bmu\, g''(\br)\, \br_t\, \xi 
    - \bmu\, g'(\br)\, \xi_t 
  \bigr) \,\eta \dx\ds
  \non
  \\
  && \leq \frac 14 \intQt |\xi_t|^2 \dx\ds
  + c \intQt ( |\xi|^2 + |\eta|^2 ) \dx\ds \,.
  \non
\Esist
Finally, \pier{recalling~\eqref{B5}, it is clear that}
\Bsist
  && H_3 := \intQt \bigl(
    ( \juerg{\xi} + \bmu\, g''(\br) - F''(\br) ) \,\xi
    - DB[\br](\xi)
    + g'(\br) \, \eta
  \bigr) \,\xi_t \dx\ds
  \non
  \\
  && \leq \frac 14 \intQt |\xi_t|^2 \dx\ds
  + c \intQt ( |\xi|^2 + |\eta|^2 ) \dx\ds \,.
\Esist
Therefore, it suffices to \pier{collect these inequalities and apply Gronwall's lemma in order}
to conclude that $\xi=0$ and $\eta=0$.

\medskip

The existence of a solution is proved in several steps.
First, we introduce an approximating problem depending on the parameter $\eps\in(0,1)$.
Then, we show well-posedness for \juerg{this} problem
and perform suitable a priori estimates.
Finally, we construct a solution to problem \accorpa{ls1}{ls4} by letting $\eps$ tend to zero.
For the sake of simplicity, in performing the uniform a~priori estimates,
we denote by $c>0$ different constants that may depend on the data
of the system but not on $\eps\in(0,1)$; the actual value of $c$ may change within formulas 
and lines.
On the contrary, the symbol $\ceps$ stands for (different) constants that can depend also on~$\eps$.
In particular, $\ceps$~is independent of the parameter~$\delta$
that enters an auxiliary problem we introduce later~on.

\step 1
We approximate $\br$ and $\bmu$
by suitable $\rhoeps,\mueps\in C^\infty(\overline Q)$ as specified below.
For every $\eps\in(0,1)$, it holds that
\Beq
  \rho_{**} \leq \rhoeps \leq \rho^{**}
  \ \hbox{in $\overline Q$}
  \aand
  \norma{\rhoeps_t}_{\LQ\infty} + \norma\mueps_{\H1{\Lx3}\cap\LQ\infty} \leq C^*,
  \label{hpapprox}
\Eeq
for some constants $\rho_{**},\rho^{**}\in(0,1)$ and $C^*>0$;
as $\eps\searrow0$, we have
\Bsist
  && \hskip -3em \rhoeps \to \br, \ \
  \rhoeps_t \to \br_t, \ \
  \mueps \to \bmu, \ \
  \hbox{in $\LQ p$, for every $p<+\infty$ and a.\,e.\ in $Q$,}
  \non
  \\
  && \aand
  \mueps_t \to \bmu_t \ \ 
  \hbox{in $\L2{\Lx3}$}.
  \label{conveps}
\Esist
\pier{In order to construct regularizing families as above, we can use, for instance, extension outside $Q$ and convolution with mollifiers.}

Next, we introduce the approximating problem of finding $(\xieps,\etaeps)$ satisfying
\Bsist
  && \xieps_t + F''(\br)\,\xieps + DB[\br](\xieps)
  = \bmu \,g''(\br)\,\xieps + g'(\br)\,\etaeps
  \quad\mbox{a.\,e.\ in }\,Q,
  \label{eqxieps}
  \\[1mm]
  && (1+2g(\rhoeps))\, \etaeps_t + g'(\rhoeps)\,\rhoeps_t\,\etaeps 
  \nonumber
  \\ 
  && \,\,
  +\, 2g'(\br)\,\mueps_t\,\xieps + \bmu g''(\br)\,\br_t\,\xieps 
	+ \bmu g'(\br)\,\xieps_t
  -\Delta\etaeps = h
  \quad\mbox{a.\,e.\ in }\,Q,
  \qquad\quad
  \label{eqetaeps}
  \\[1mm]
  && \dn\etaeps=0 \quad\mbox{a.\,e.\ on }\, \Sigma,
  \label{bceps}
  \\[1mm]
  && \etaeps(0)=\xieps(0)=0 \quad\mbox{a.\,e.\ in }\,\oma.
  \label{iceps}
\Esist
In order to solve \accorpa{eqxieps}{iceps},
we introduce the spaces
\Beq
  \calV := H\times V
  \aand
  \calH := H\times H,
  \non
\Eeq
and present our problem in the form
\Beq
  \frac d{dt} \, (\xi,\eta) + \Aeps (\xi,\eta) = f
  \aand
  (\xi,\eta)(0) = (0,0),
  \non
\Eeq
in the framework of the Hilbert triplet $(\calV,\calH,\calV')$.
We look for a weak solution
and aim at applying \cite[Thm.~3.2, p.~256]{Baiocchi}.
To this end, we have to split $\,\Aeps\,$ in the form $\,\pier{\calQ^\eps}+\Reps$,
where $\pier{\calQ^\eps}$ is the uniformly elliptic principal part
and the remainder $\Reps$ satisfies the requirements
\cite[(4.4)--(4.5), p.~259]{Baiocchi}.
We notice at once that \juerg{these} conditions are trivially fulfilled whenever
\Bsist
  && \hskip-3em
  \Reps = (\Reps_1 , \Reps_2) \in \calL (\L2\calH,\L2\calH)\,,
  \label{baioA}
  \\
  && \hskip-3em
  \left|
    \intQt \bigl( \Reps_1(v,w) \, v + \Reps_2(v,w) \, w \bigr) \dx\ds
  \right|
  \leq C_{\Reps} \intQt \bigl( |v|^2 + |w|^2 \bigr) \ds\,,
  \qquad
  \label{baioB}
\Esist
for some constant $C_{\Reps}$, and every $v,w\in\L2H$ and $t\in[0,T]$.
In order to present \accorpa{eqxieps}{iceps} in the desired form,
we multiply \eqref{eqetaeps} by $\aeps:=1/(1+\juerg{2}g(\rhoeps))$ and notice that
\Beq
  - \aeps \Delta\etaeps = - \div(\aeps\nabla\etaeps) + \nabla\aeps \cdot \nabla\etaeps .
  \non
\Eeq
As $\aeps\geq\alpha:=1/(1+2\sup g)$ and $\nabla\aeps$ is bounded,
we can fix a real number $\lambdaeps>0$ such~that
\Beq
  \iO \bigl(
    \aeps(t) |\nabla w|^2 + (\nabla\aeps(t) \cdot \nabla w) w + \lambdaeps |w|^2
  \bigr)\juerg{\dx} \,\geq\, \frac \alpha 2 \, \normaV w^2
  \label{coerc}
\Eeq
for every $w\in V$ and $t\in[0,T]$.
Next, we replace $\xieps_t$ in \eqref{eqetaeps} by using \eqref{eqxieps}.
Therefore, we see that a possible weak formulation of \accorpa{eqxieps}{bceps} is given~by
\Bsist
  &&  \pier{\iO \xieps_t(t) v \dx  + {}_{V'}\< \etaeps_t(t) , w >_V} 
  + {}_{\calV'} \< \calQ^\eps(t) (\xieps,\etaeps)(t) , (v,w) >_{\calV}
  \non
  \\
  && \quad \hbox{}
  + \iO \bigl( \Reps_1(\xieps,\etaeps)\juerg{(t)} v + 
	\Reps_2(\xieps,\etaeps)\juerg{(t)}\, w \bigr) \dx
  = \iO \aeps\juerg{(t)}\, h\juerg{(t)}\, w\juerg{\dx}
  \non
  \\
  && \hbox{for a.\ a. $t\in(0,T)$ and every $(v,w)\in\calV$,}
  \label{lseps}
\Esist
where the symbols $\<\cdot,\cdot>$ stand for the duality pairings
and $\calQ^\eps$ and $\Reps_i$ have the meaning explained below.
The time-dependent operator $\calQ^\eps(t)$ from $\calV$ into $\calV'$ is defined~by
\Bsist
  && {}_{\calV'} \< \calQ^\eps(t) (\hat v,\hat w) , (v,w) >_{\calV}
  \non
  \\
  && = \iO \bigl( \hat v\, v + \aeps(t)\, \nabla\hat w \cdot \nabla w + (\nabla\aeps(t)\cdot\nabla\hat w)\, w + \lambdaeps\,\hat w\, w \bigr) \dx
  \quad
  \label{formaeps}
\Esist
for every $(\hat v,\hat w),(v,w)\in\calV$ and $t\in[0,T]$.
By construction, the bilinear form given by the \rhs\ of \eqref{formaeps}
is continuous on $\calV\times\calV$, depends smoothly on time, 
and is $\calV$-coercive uniformly with respect to~$t$ (see~\eqref{coerc}).
The operators 
$$\pier{\Reps_i\in\calL(\L2\calH,\L2H)}$$ 
account for the term $\lambdaeps\etaeps$ 
that has to be added also to the \rhs\ of~\eqref{eqetaeps}
and for all the contributions to the equations that have not been considered in the principal part.
They have the form
\Beq
  (\Reps_i(v,w))(t)
  = \aeps_{i1}(t) v + \aeps_{i2}(t) w + \aeps_{i3}(t) \bigl( DB[\br](v) \bigr)(t)
  \label{restieps}
\Eeq
for $(v,w)\in \L2\calH$, with some coefficients $\aeps_{ij}\in\LQ\infty$.
Therefore, \pier{by virtue of}~\eqref{B5}, we see that
\Bsist
  && \intQt \bigl( \Reps_1(v,w) \, v + \Reps_1(v,w) \, w \bigr) \dx\ds
  \non
  \\
  && \leq c \intQt (|v|^2+|w|^2) \dx\ds + c \,\norma{DB[\br](v)}_{L^2(Q_t)}^2
  \non
  \\
  && \leq c \intQt (|v|^2+|w|^2) \dx\ds\,,
  \non
\Esist
for every $(v,w)\in\L2\calH$ and every $t\in[0,T]$.
Thus, the conditions \accorpa{baioA}{baioB} are fulfilled,
and the result of \cite{Baiocchi} mentioned above can be applied.
We conclude that the Cauchy problem for \eqref{lseps} has a unique solution $(\xieps,\etaeps)$ satisfying
\Bsist
  && (\xieps,\etaeps) \in \H1{\calV'} \cap \L2\calV,
  \quad \hbox{i.\,e.,}
  \non
  \\
  && \xieps \in \H1H
  \aand
  \etaeps \in \H1{V'} \cap \L2V.
  \non
\Esist
On the other hand, \juerg{this} solution has to satisfy 
\Beq
  \< \pt\etaeps , w > + \iO \aeps\,\nabla\etaeps \cdot \nabla w \dx = \iO \phi_\eps\, w \dx
  \quad \hbox{a.\,e.\ in (0,T), for every $w\in V$,}
  \non
\Eeq
for some $\phi_\eps\in\LQ2$. From standard elliptic regularity, it follows 
that $\etaeps\in\H1H\cap\L2W$.

In the next steps,  
besides of Young's inequality, we make repeated use of  
the global estimates \pier{\eqref{ssbounds1}, \eqref{ssbounds2},} and \eqref{ssbounds3}, for $\br$ and~$\bmu$, 
without further reference.

\step 2
For convenience, we refer to \juerg{Eqs.} \accorpa{eqxieps}{bceps} 
(\juerg{using} the language that is proper for strong solutions),
but it is understood that they are meant in the variational sense~\eqref{lseps}.
We add $\xieps$ and $\etaeps$ to both sides of \eqref{eqxieps} and~\eqref{eqetaeps}, respectively; then, we multiply the resulting equalities by $\xieps_t$ and~$\etaeps$,
integrate over $Q_t$, and sum~up.
By observing that
\Beq
  (1+2g(\rhoeps))\, \etaeps_t\, \etaeps + g'(\rhoeps)\, \rhoeps_t \,|\etaeps|^2
  = \pt \bigl[\juerg{\bigl( {\smfrac 12} + g(\rhoeps)\bigr)}\, |\etaeps|^2 \bigr],
  \non
\Eeq
and recalling that $g\geq0$, we obtain
\Beq
  \frac 12 \iO |\xieps(t)|^2 \dx
  + \intQt |\xieps_t|^2 \dx\ds
  + \frac 12 \iO |\etaeps(t)|^2 \dx
  + \iot \normaV{\etaeps(s)}^2 \ds
  \leq \sum_{j=1}^3 I_j,
  \non
\Eeq
where the terms $I_j$ are defined and estimated as follows. 
\pier{In view of~\eqref{B5},} we \juerg{first} \pier{infer} \juerg{that}
\begin{align}
  I_1 :&= \intQt \bigl(
  \xieps 
  - F''(\br)\, \xieps 
  - DB[\br](\xieps) 
  + \bmu \,g''(\br) \,\xieps 
  + g'(\br) \,\etaeps
  \bigr) \,\xieps_t \dx\ds
  \non
  \\
  & \leq \pier{\,\frac 14 \intQt |\xieps_t|^2 \dx\ds
  + c \intQt \bigl( |\xieps|^2 + |\etaeps|^2 \bigr) \dx\ds \, .} 
  \non
\end{align}
Next, we have
\begin{align}
  I_2 :&= \intQt \bigl(
  \etaeps
  - \bmu\, g''(\br)\, \br_t\, \xieps
  - \bmu g'(\br)\, \xieps_t 
  + h
  \bigr) \,\etaeps \dx\ds
  \non
  \\
  & \leq \,\frac 14 \intQt |\xieps_t|^2 \dx\ds
  + c \intQt \bigl( |\xieps|^2 + |\etaeps|^2 \bigr) \dx\ds 
  + c \,.
  \non
\end{align}
Finally, by virtue of the H\"older and Sobolev inequalities, we have
\begin{align}
  I_3 :&= - \intQt 2 g'(\br)\, \mueps_t\, \xieps\, \etaeps\juerg{\dx\ds}
  \leq c \iot \norma{\mueps_t(s)}_3 \, \norma{\xieps(s)}_2 \, \norma{\etaeps(s)}_6 \ds
  \non
  \\
  & \leq \frac 12 \iot \normaV{\etaeps(s)}^2 \ds
  + c \iot \norma{\mueps_t(s)}_3^2 \, \norma{\xieps(s)}_2^2 \ds \,.
  \non
\end{align}
At this point, we recall all the inequalities we have proved,
notice that \eqref{hpapprox} implies that the function
$s\mapsto\norma{\mueps_t(s)}_3^2$ is bounded in $L^1(0,T)$, and apply the Gronwall lemma.
We obtain
\Beq
  \norma\xieps_{\H1H} + \norma\etaeps_{\L\infty H\cap\L2V} \leq c \,.
  \label{primastima}
\Eeq

\step 3
We would \juerg{now like to test} \eqref{eqxieps} by $(\xieps)^5$.
However, \juerg{this} function is not admissible, unfortunately.
Therefore, we introduce a suitable approximation.
To start with, we consider the Cauchy problem
\Beq
  \hat\xi_t + b\, \hat\xi + \juerg{L( \hat\xi)} = \feps
  \aand
  \hat\xi(0) = 0,
  \label{hatxi}
\Eeq
where we have set, for brevity,
\Beq
  b := F''(\br) - \bmu\, g''(\br), \quad
  L := DB[\br],
  \aand
  \feps := g'(\br)\, \etaeps .
  \label{approxdelta}
\Eeq
By \eqref{eqxieps}, $\hat\xi:=\xieps$ is a solution belonging to $\H1H$.
On the other hand, such a solution is unique.
Indeed, multiplying by~$\hat\xi$ the corresponding homogeneous equation
(i.\,e., $\feps$~is replaced by~$0$),
and invoking \eqref{B5} and Gronwall's lemma, one immediately obtains that $\hat\xi=0$.
We conclude that $\hat\xi:=\xieps$ is the unique solution to~\eqref{hatxi}.

At this point, we approximate $\xieps$ by the solution to a problem 
depending on the parameter $\delta\in(0,1)$, in addition.
Namely, we look for $\xied$ satisfying the parabolic-like equation
\Beq
  \xied_t - \delta \,\Delta\xied + \bdelta\, \xied + \juerg{L(\xied)} = \feps,
  \label{pbled}
\Eeq
complemented with the Neumann boundary condition $\dn\xied=0$ and the initial condition $\xied(0)=0$.
In \eqref{pbled}, $\bdelta$~is an approximation of $b$ belonging to $C^\infty(\overline Q)$ \juerg{that satisfies}
\Beq
  \norma\bdelta_{\LQ\infty} \leq c,
  \aand
  \bdelta \to b 
  \quad \hbox{a.\,e.\ in $Q$ as $\delta\searrow0$}.
  \label{convdelta}
\Eeq
\juerg{This} problem has a unique weak solution $\xied\in\H1{V'}\cap\L2V$,
as one easily sees by arguing as we did for 
\juerg{the} more complicated system \accorpa{eqxieps}{iceps}
and applying \cite[Thm.~3.2, p.~256]{Baiocchi}.

\juerg{We now aim to show} that $\xied$ is bounded.
To this end, we introduce the operator $\Adelta\in\calL(V,V')$
defined by 
$$\<\Adelta v,w>:=\iO(\delta\,\nabla v\cdot\nabla w+v\,w)\dx \quad\mbox{for every $v,w\in V$,}
$$
and observe that $\Adelta$ is an isomorphism.
Moreover, \juerg{Eq.} \eqref{pbled}, complemented with the boundary and initial conditions,
can be written~as
\Beq
  \xied_t + \Adelta\xied = \fed := \feps - (1+\bdelta)\, \xied + \juerg{L(\xied)} 
  \aand 
  \xied(0) = 0 \,.
  \label{pbledbis}
\Eeq
Now, by also accounting for~\eqref{B6},
we notice that $\feps$, $\xied$, $\bdelta\xied$,
and $\juerg{L(\xied)}$, all belong to $\L2V$.
Hence, $\fed\in\L2V$, so that $\Adelta\fed\in\L2{V'}$, 
and we can consider the unique solution $\zed\in\H1{V'}\cap\L2V$
to the problem
\Beq
  \zed_t + \Adelta\zed = \Adelta\fed
  \aand
  \zed(0) = 0 \,.
  \non
\Eeq
Now, $\Adelta^{-1}\zed$ satisfies 
\Beq
  (\Adelta^{-1}\zed)_t + \Adelta (\Adelta^{-1}\zed) = \Adelta^{-1} \Adelta\fed 
	= \juerg{\fed}
  \aand
  (\Adelta^{-1}\zed)(0) = 0,
  \non
\Eeq
so that a \juerg{comparison} with \eqref{pbledbis} shows that $\xied=\Adelta^{-1}\zed$,
by uniqueness.
Since $\zed\in\L\infty H$, and $\Adelta^{-1}(H)=W$ by elliptic regularity,
we deduce that $\xied\in\L\infty W$. 
Therefore, $\xied$~is bounded, \juerg{as claimed}.

\juerg{Consequently,} $(\xied)^5$ is an admissible test function,
since it \gian{clearly} belongs to \pier{the space} $\L2V$.
By multiplying \eqref{pbled} by $(\xied)^5$ and integrating over~$Q_t$, we obtain that
\Beq
  \frac 16 \iO |\xied(t)|^6 \juerg{\dx}\,+\, 5\,\delta \intQt 
	|\xied|^4\, |\nabla\xied|^2 \dx\ds
  = \sum_{j=1}^3 K_j,
  \non
\Eeq
where the terms $K_j$ are defined and estimated as follows.
First, recalling \juerg{\eqref{convdelta}}, we deduce that
\Beq
  K_1 := - \intQt \bdelta\,\xied\, (\xied)^5 \dx\ds
  \,\leq\, c \intQt |\xieps|^6 \dx\ds \,.
  \non
\Eeq
On the other hand, H\"older's inequality, and assumption \eqref{B5} with $p=6$, imply that
\begin{align}
  K_2 :&= - \intQt \juerg{L(\xied)}\, (\xied)^5 \dx\ds
  \,\leq\, c \, \norma{L\xied}_{L^6(Q_t)} \, \norma{(\xied)^5}_{L^{6/5}(Q_t)}
  \non
  \\
  & \leq c \, \norma\xied_{L^6(Q_t)} \, \norma\xied_{L^6(Q_t)}^5
  \,=\, c \intQt |\xied|^6 \dx\ds. 
  \non
\end{align}
Finally, also invoking Sobolev's inequality, we see that
\begin{align}
  K_3 :&= \intQt \feps\, (\xied)^5 \dx\ds
  \,\leq\, c \iot \norma{\etaeps(s)}_6 \, \norma{(\xied(s))^5}_{6/5} \ds
  \non
  \\
  & \pier{{}\leq} \,c \iot \norma{\etaeps(s)}_6 \, \norma{\xied(s)}_6^5 \ds
  \leq c \iot \normaV{\etaeps(s)} \, \bigl( 1 + \norma{\xied(s)}_6^6 \bigr) \ds \,.
  \non
\end{align}
\juerg{Collecting the above estimates}, and noting that
the function $s\mapsto\normaV{\etaeps(s)}$ is bounded in $L^1(0,T)$
by \eqref{primastima},
we can apply the Gronwall lemma to conclude that
\Beq
  \norma\xied_{\L\infty\Lsei} \,\leq\, c \,.
  \label{stimaed}
\Eeq
At this point, we quickly show that $\xied$ converges to~$\xieps$ as $\delta\searrow0$,
at least for a subsequence.
Indeed, one multiplies \eqref{pbled} first by $\xied$, and then by~$\xied_t$,
and proves that 
\Beq
  \norma\xied_{\H1H\cap\L\infty V} \,\leq\, \ceps,
  \non
\Eeq
uniformly with respect to~$\delta$.
Then, by weak compactness and \eqref{convdelta}
(which implies convergence \pier{of $b^\delta $ to $b$} in $L^p(Q)$ for every $p<+\infty$), 
it~is \sfw\ to see that
$\xied$ converges to a solution $\hat\xi$ to the problem associated with~\eqref{hatxi}.
As $\hat\xi=\xieps$, we have proved what we have claimed.
This, and \eqref{stimaed}, yield that
\Beq
  \norma\xieps_{\L\infty\Lsei} \,\leq\, c \,.
  \label{secondastima}
\Eeq

\step 4
At this point, we can multiply \eqref{eqetaeps} by $\etaeps_t$ and integrate over~$Q_t$.
\pier{By recalling that $g\geq0$, we obtain} 
\Beq
  \intQt |\etaeps_t|^2 \dx\ds
  + \frac 12 \iO |\nabla\etaeps(t)|^2 \dx
  \,\leq\, \sum_{j=1}^3 L_j,
  \non
\Eeq
where each term $L_j$ is defined and estimated below. First,
by taking advantage of \eqref{secondastima} and \eqref{hpapprox} for~$\mueps_t$, we have
\begin{align}
  L_1 :&= - \intQt 2 g'(\br)\, \mueps_t \,\xieps\, \etaeps_t \dx\ds
  \leq c \iot \norma{\mueps_t(s)}_3 \, \norma{\xieps(s)}_6 \, \norma{\etaeps_t(s)}_2 \ds
  \non
  \\
  & \leq\, c \iot \norma{\mueps_t(s)}_3 \, \norma{\etaeps_t(s)}_2 \ds
  \,\leq \,\frac 14 \intQt |\etaeps_t|^2 \dx\ds
  + c \ioT \norma{\mueps_t(s)}_3^2 \ds
  \non
  \\
  & \leq\, \frac 14 \intQt |\etaeps_t|^2 \dx\ds + c \,.
  \non
\end{align}
Next, using \eqref{hpapprox} for $\rhoeps_t$ \pier{and \eqref{primastima}}, we obtain that
\begin{align}
  L_2 := - \intQt g'(\rhoeps) \rhoeps_t \etaeps \etaeps_t \dx\ds
  \leq &\ \frac 14 \intQt |\etaeps_t|^2 \dx\ds
  + c \intQt |\etaeps|^2 \dx\ds \non\\
  \pier{\leq} &\ \pier{\frac 14 \intQt |\etaeps_t|^2 \dx\ds
  + c  }
  \,.
  \non
\end{align}
Finally, in view of \eqref{primastima}, we have
\begin{align}
  L_3 :&= \intQt \bigl( - \bmu\, g''(\br)\,\pier{\br_t}\, \xieps - 
	\bmu\, g'(\br)\, \xieps_t + h \bigr)\, \etaeps_t \dx\ds
  \non
  \\
  & \leq\, \frac 14 \intQt |\etaeps_t|^2 \dx\ds
  + c \intQt \bigl( |\xieps|^2 + |\xieps_t|^2 + 1 \bigr) \dx\ds
  \non
  \\
  & \leq\, \frac 14 \intQt |\etaeps_t|^2 \dx\ds + c \,.
  \non
\end{align}
By collecting \juerg{the above} estimates, we conclude that
\Beq
  \norma{\etaeps_t}_{\L2H} + \norma\etaeps_{\L2V}
  \,\leq\, c \,.
  \label{terzastima}
\Eeq
As a consequence, we can estimate $\Delta\etaeps$ in $\LQ2$, just by 
comparison in~\eqref{eqetaeps}.
Using standard elliptic regularity, we deduce that
\Beq
  \norma\etaeps_{\L2W} \,\leq\, c \,.
  \label{quartastima}
\Eeq

\step 5
At this point, we are ready to prove the existence part of the statement.
Indeed, the estimates \eqref{primastima} and \accorpa{secondastima}{quartastima}
yield that
\Bsist
  & \xieps \to \xi
  & \quad \hbox{weakly star in $\H1H\cap\L\infty\Lsei$},
  \non
  \\
  & \etaeps \to \eta
  & \quad \hbox{weakly star in $\H1H\cap\L\infty V\cap\L2W$},
  \non
\Esist
as $\eps\searrow0$, at least for a subsequence.
By accounting for \eqref{conveps} and the Lipschitz continuity of $g$ and~$g'$,
it is \sfw\ to see that $(\xi,\eta)$ is a solution to problem \accorpa{ls1}{ls4}.
This completes the proof.\qed

\vspace{5mm}
We are now prepared to show that $\cs$ is directionally differentiable. We 
have the following result:

\vspace{5mm}
{\sc Theorem 3.2:} \quad {\em Suppose that the general hypotheses} {\bf (A1)}--{\bf (A5)}
{\em are satisfied, and let $\bu\in\uad$ be given and $(\br,\bmu)=\cs(\bu)$. Moreover,
let $h\in\cx$ be \pier{a function} such that $\bu+\overline{\lambda} h\in\uad$ for some
$\overline{\lambda}>0$. Then the directional derivative $\delta\cs(\bu;h)$ of $\cs$ at $\bu$ in the
direction $h$ exists in the space $(\cy, \|\cdot\|_\cy)$, and we have
$\delta\cs(\bu;h)=(\xi,\eta)$, where $(\xi,\eta)$ is the unique solution to the linearized
system} \eqref{ls1}--\eqref{ls4}. 

\vspace{3mm}
{\sc Proof:} We have 
$\bu+\lambda h\in\uad$ for $0<\lambda\le\blam$, since $\uad$ is convex. We put, for any such $\lambda$,
$$\ula:=\bu+\lambda h, \quad (\rla,\mula):=\cs(\ula), \quad
\yla:=\rla-\br-\lambda\xi, \quad \zla:=\mula-\bmu-\lambda\eta.$$
Notice that
$\,(\rla,\mula)$ and $(\br,\bmu)$ fulfill the global bounds \eqref{ssbounds1}, \eqref{ssbounds2}, and 
\eqref{ssbounds3}, and that \pier{$(\yla, \zla) \in \cal Y$ for all $\lambda\in[0,\blam]$}.
Moreover, by virtue of Theorem 2.2, we have the estimate
\begin{align}
\label{stabu2}
&\|\rla-\br\|_{H^1(0,t;H)\cap L^\infty(0,t;L^6(\oma))}
\,+\,\|\mula-\bmu\|_{H^1(0,t;H)\cap L^\infty(0,t;V)\cap L^2(0,t;\pier{W})}\nonumber\\[2mm]
&\le\,K_2^*\,\lambda\|h\|_{\lzht}\,.
\end{align} 
We also notice that, owing to \eqref{ssbounds1} and the assumptions on $F$ and $g$, it follows from
Taylor's theorem that
\begin{align}
\label{taylor1}
&\left|F'(\rla)-F'(\br)-\lambda F''(\br)\xi\right|\,\le\,c\left|\yla\right|\,+\,c\left|\rla-\br\right|^2\quad
\mbox{a.\,e. in }\,Q,\\[1mm]
\label{taylor2}
&\left|g(\rla)-g(\br)-\lambda g'(\br)\xi\right|\,\le\,c\left|\yla\right|\,+\,c\left|\rla-\br\right|^2\quad
\mbox{a.\,e. in }\,Q,\\[1mm]
\label{taylor3}
&\left|g'(\rla)-g'(\br)-\lambda g''(\br)\xi\right|\,\le\,c\left|\yla\right|\,+\,c\left|\rla-\br\right|^2\quad
\mbox{a.\,e. in }\,Q,
\end{align}
where, here and in the remainder of the proof, we denote by $c$ constants that may depend on the data
of the system but not on $\lambda\in [0,\blam]$; the actual value of $c$ may change within formulas 
and lines. Moreover, by the Fr\'echet differentiability of~$B$ 
(recall assumption {\bf (A4)(vi)} \gianni{and the fact that, for $\overline v,v\in L^2(Q)$,
the restrictions of $B[v]$ and $DB[\overline v](v)$ to $Q_t$
depend only on~$v|_{Q_t}$}), 
we have \pier{(cf.~\eqref{stabu2})}
\beq
\label{taylor4}
\|B[\rla]-B[\br]-\lambda \,DB[\br](\xi)\|_{L^2(\gianni{Q_t})}\,
  \le\, \gianni{c\,\|y^\lambda\|_{L^2(Q_t)}+ R\left(\lambda\,\|h\|_{L^2(Q_t)}\right)},
\eeq
with a function $\gianni R:(0,+\infty)\to (0,+\infty)$ satisfying \,\,$\lim_{\sigma\searrow 0}\gianni R(\sigma)/\sigma=0$.
\gianni{As we want to prove that $\delta S(\bu;h)=(\xi,\eta)$,
according to the definition of directional differentiability, 
we need to show that}
\begin{align}
\label{Frechet}
0\,&=\,\lim_{\lambda\searrow 0}\,\frac {\|\cs(\bu+\lambda h)-\cs(\bu)-\lambda \,\gianni{(\xi,\eta)}\|_{\cy}}
\lambda\nonumber\\[1mm]
&=\,\lim_{\lambda\searrow 0}\,\frac{\left\|\yla\right\|_{H^1(0,T;H)}\,
+\,\left\|\zla\right\|_{\pier{L^\infty(0,T;H)\cap L^2(0,T;V)}}} \lambda \,\,.
\end{align}

To begin with, using the state system \eqref{ss1}--\eqref{ss4} and the linearized system \eqref{ls1}--\eqref{ls4},
we easily verify that for $0<\lambda\le\blam$ the pair $(\zla,\yla)$ is a strong solution
to the system
\begin{align}
\label{yz1}
&(1+2g(\br))\zla_t\,+\,g'(\br)\br_t\zla\,+\,\bmu\,g'(\br)\yla_t-\Delta \zla\nonumber\\
&=-\,2\left(g(\rla)-g(\br)\right)\left(\mula_t-\bmu_t\right)\,
-\,2\,\bmu_t\left(g(\rla)-g(\br)-
\lambda g'(\br)\xi\right)\nonumber\\
&\quad \,-\,\bmu\,\br_t\left(g'(\rla)-g'(\br)-\lambda g''(\br)\xi\right)
\,-\,\bmu\left(g'(\rla)-g'(\br)\right)\left(\rla_t-\br_t\right)\nonumber\\
&\quad\,-\left(\mula-\bmu\right)\left[\bigl(g'(\rla)-g'(\br)\bigr)\,\br_t
+ g'(\rla)\left(\rla_t-\br_t\right)\right] \quad\mbox{a.\,e. in }\, Q,\\[1mm]
\label{yz2}
&\yla_t=-\left(F'(\rla)-F'(\br)-\lambda F''(\br)\xi\right)-\left(B[\rla]-B[\br]-\lambda DB[\br](\xi)\right)
\nonumber\\
&\hspace*{9mm} +\,g'(\br)\zla\,+\,\bmu\left(g'(\rla)-g'(\br)-\lambda g''(\br)\xi\right)\nonumber\\
&\hspace*{9mm} +\left(\mula-\bmu\right)\left(g'(\rla)-g'(\br)\right)  \quad\mbox{a.\,e. in}\,Q,\\[1mm]
\label{yz3}
&\hspace*{3cm}\dn\zla=0 \quad\mbox{a.\,e. on }\,\Sigma,\\[1mm]
&\hspace*{2.9cm}\zla(0)=\yla(0)=0 \quad\mbox{a.\,e. in }\,\Omega.
\end{align}

In the following, we make repeated use of the mean value theorem and 
of the global estimates \eqref{ssbounds1}, \eqref{ssbounds2},
\eqref{ssbounds3}, and \pier{\eqref{stabu2}}, without further reference. 
\gianni{For the sake} of a better readability, we will
omit the superscript $\lambda$ of the quantities $\yla$, $\zla$ during the estimations, writing it only
at the end of the respective estimates.

\vspace{3mm}
\underline{\sc Step 1:}  \,\,\,Let $t\in (0,T]$ be fixed. First, observe that
$$\pt\Bigl(\Bigl(\frac 12 + g(\br)\Bigr)z^2\Bigr)\,=
\,(1+2\,g(\br))\,z\,z_t \,+\,g'(\br)\,\br_t\,z^2\,.
$$ 
Hence, multiplication of \eqref{yz1} by $z$ and  integration over $Q_t$ yields the estimate
\beq
\label{p321}
\ioma \Bigl(\frac 12 + g(\br(t))\Bigr)z^2(t)\dx\,+\,\txinto|\nabla z|^2\dx\ds\,\le\,c\sum_{j=1}^7 |I_j|\,,
\eeq
where the quantities $I_j$, $1\le j\le 7$, are specified and estimated as follows: at first, Young's
inequality shows that, for every $\gamma>0$ (to be chosen later),
\beq
\label{p322}
I_1:=\,-\txinto \bmu\,g'(\br)\,y_t\,z\dx\ds\,\le\,\gamma\txinto y_t^2\dx\ds\,+\,\frac c\gamma\txinto z^2\dx\ds\,.
\eeq  
Moreover, we have, by H\"older's and Young's inequalities and \pier{\eqref{stabu2}},
\begin{align}
\label{p323}
I_2:&=\,-2\txinto\left(g(\rla)-g(\br)\right)\left(\mula_t-\bmu_t\right)z\dx\ds\nonumber\\[1mm]
&\le\,c\tint\|\rla(s)-\br(s)\|_6\,\|\mula_t(s)-\bmu_t(s)\|_2\,\|z(s)\|_3\ds\nonumber\\[1mm]
&\le\,
c\,\|\rla-\br\|_{L^\infty(0,t;L^6(\oma))}\,\|\mula-\bmu\|_{H^1(0,t;H)}\,\|z\|_{\lzvt}\nonumber\\[1mm]
&\le\,\gamma\,\|z\|_{\lzvt}^2\,+\,\frac c\gamma \,\lambda^4\,.
\end{align}
Next, we employ \eqref{taylor2}, the H\"older and Young inequalities, and \pier{\eqref{stabu2}}, to infer that
\begin{align}
\label{p324}
I_3:&=\,-2\txinto\bmu_t\left(g(\rla)-g(\br)-\lambda g'(\br)\xi\right)z\dx\ds\nonumber\\[1mm]
&\le \,c\txinto|\bmu_t|\left(|y|\,+\,|\rla-\br|^2\right)\,|z|\dx\ds\nonumber\\[1mm]
&\le \,c\tint\|\bmu_t(s)\|_6\left(\|y(s)\|_2\,\|z(s)\|_3\,+\,\|\rla(s)-\br(s)\|_6^2
\,\|z(s)\|_2\right)\ds\nonumber\\[1mm]
&\le\,\gamma\tint\|z(s)\|_V^2\ds\,+\,\frac c\gamma\tint\|\bmu_t(s)\|_V^2\,\|y(s)\|_H^2\ds\nonumber\\[1mm]
&\qquad +\,c\tint\|\bmu_t(s)\|_V^2\,\|z(s)\|_H^2\ds\,+\,c\,\|\rla-\br\|_{L^\infty(0,t;V)}^4\nonumber\\[1mm]
&\le \,\gamma\tint\|z(s)\|_V^2\ds\,+\,\bigg(1+\frac c\gamma\bigg)\!\!\tint\|\bmu_t(s)\|_V^2
\bigl
(\|y(s)\|_H^2\,+\,\|z(s)\|^2_H\bigr)\ds \pier{{}+c\,\lambda^4,}
\end{align}
where we observe that, in view of \eqref{ssbounds2}, the mapping 
$\,s\mapsto \|\bmu_t(s)\|_V^2\,$ belongs to
$L^1(0,T)$. 
Likewise, utilizing \eqref{ssbounds2}, \eqref{taylor3}, \pier{\eqref{stabu2}}, and the H\"older and Young inequalities, \pier{it is straightforward to deduce} that
\begin{align}
\label{p325}
I_4:&=\,-\txinto\bmu\,\br_t \left(g'(\rla)-g'(\br)-\lambda g''(\br)\xi\right)z
\dx\ds\nonumber\\[1mm]
&\le\,c\txinto\left(|y|+|\rla-\br|^2\right)|z|\dx\ds\nonumber\\[1mm]
&\le\,c\txinto\left(y^2+z^2\right)dx\ds\,+\,c\tint\|\rla(s)-\br(s)\|_4^2
\,\|z(s)\|_2\ds\nonumber\\[1mm]
&\le \,c\txinto\left(y^2+z^2\right)dx\ds\,+\,c\,\lambda^4\,.
\end{align}

In addition, arguing similarly, we have
\begin{align}
\label{p326}
I_5:&=\,-\txinto\bmu\left(g'(\rla)-g'(\br)\right)\left(\rla_t-\br_t\right)z\dx\ds
\nonumber\\[1mm]
&\le\,c\tint\|\rla(s)-\br(s)\|_6\,\|\rla_t(s)-\br_t(s)\|_2\,\|z(s)\|_3\ds\nonumber\\[1mm]
&\le\,c\,\|\rla-\br\|_{L^\infty(0,t;L^6(\oma))}
\,\|\rla-\br\|_{\pier{H^1(0,t;H)}}\,\|z\|_{\lzvt}\nonumber\\[1mm]
&\le\,\gamma\tint\|z(s)\|_V^2\ds\,+\,\frac c\gamma\,\lambda^4\,,
\end{align}
as well as
\begin{align}
\label{p327}
I_6:&=\,-\txinto\br_t\left(\mula-\bmu\right)\left(g'(\rla)-g'(\br)\right)z\dx\ds\nonumber\\[1mm]
&\le\,c\txinto\left|\mula-\bmu\right|\left|\rla-\br\right||z|\dx\ds\nonumber\\[1mm]
&\le\,c\tint\|\rla(s)-\br(s)\|_6\,\|\mula(s)-\bmu(s)\|_3\,\|z(s)\|_2\ds\nonumber\\[1mm]
&\le\,c\txinto z^2\dx\ds\,+\,c\,\lambda^4\,.
\end{align}

Finally, we find that
\begin{align}
\label{p328}
I_7:&=\,-\txinto \left(\mula-\bmu\right)g'(\rla)\left(\rla_t-\br_t\right)z\dx\ds\nonumber\\[1mm]
&\le\,c\tint\|\mula(s)-\bmu(s)\|_6\,\|\rla_t(s)-\br_t(s)\|_2\,\|z(s)\|_3\ds\nonumber\\[1mm]
&\le\,c\,\|\mula-\bmu\|_{L^\infty(0,t;V)}\,\|\rla-\br\|_{H^1(0,t;H)}\,\|z\|_{\lzvt}\nonumber\\[1mm]
&\le\,\gamma\tint\|z(s)\|_V^2\ds\,+\,\frac c\gamma\,\lambda^4\,.
\end{align}

In conclusion, combining the estimates \eqref{p321}--\eqref{p328}, 
and choosing $\gamma=\frac 18$,
we have shown that
\begin{align}
\label{p329}
&\frac 12\,\left\|\zla(t)\right\|_H^2\,+\,\frac 12\tint\left\|\zla(s)\right\|^2_V ds\,
\le\,\frac 18\txinto\left|\yla_t\right|^2 dx\ds\,+\,c\,\lambda^4\nonumber\\[1mm]
&\hspace*{1cm} +\,c\tint\left(1+\left\|\bmu_t(s)\right\|_V^2\right)
\left(\left\|\yla(s)\right\|_H^2
\,+\left\|\zla(s)\right\|_H^2\right)ds\,.
\end{align}

\vspace{3mm}
\underline{\sc Step 2:} \,\,\,Let $t\in (0,T]$ be fixed. We add $y$ to 
both sides of \eqref{yz2},
multiply the resulting identity by $y_t$, and integrate over $Q_t$ to obtain
\beq
\label{p3210}
\txinto y_t^2\dx\ds\,+\,\frac 12 \|y(t)\|_H^2\,\le\,\sum_{j=1}^6|J_j|\,, 
\eeq  
where the terms $J_j$, $1\le j\le 6$, are specified and estimated as follows: at first, 
we have, for every $\gamma>0$ (to be specified later),
\beq
\label{p3211}
J_1:=\txinto y\,y_t\dx\ds\,\le\,\gamma\txinto y_t^2\dx\ds\,+\,\frac c\gamma\txinto y^2\dx\ds\,.
\eeq
Then, we employ \pier{\eqref{ssbounds2}, \eqref{ssbounds3}, \eqref{stabu2}, and \eqref{taylor1},} as well as H\"older's and
Young's inequalities, to obtain the estimate
\begin{align}
\label{p3212}
J_2:&=\,-\txinto \left(F'(\rla)-F'(\br)-\lambda F''(\br)\xi\right)y_t\dx\ds\nonumber\\[1mm]
&\le\,c\txinto\left( |y|+|\rla-\br|^2\right)|y_t|\dx\ds\nonumber\\[1mm]
&\le\,c\tint\left(\|y(s)\|_2\,+\,\|\rla(s)-\br(s)\|_4^2\right)\|y_t(s)\|_2\ds\nonumber\\[1mm]
&\le\,\gamma\txinto y_t^2\dx\ds\,+\,\frac c\gamma\txinto y^2\dx\ds\,+\,\frac c\gamma\,\lambda^4
\,.
\end{align}

By the same token, this time invoking \eqref{taylor3}, we find that
\begin{align}
\label{p3213}
J_3:&=\txinto \bmu\left(g'(\rla)-g'(\br)-\lambda g''(\br)\xi\right) y_t\dx\ds\nonumber\\[1mm]
&\le\,\gamma\txinto y_t^2\dx\ds\,+\,\frac c\gamma\txinto y^2\dx\ds\,+\,\frac c\gamma\,\lambda^4\,.
\end{align}
 
Moreover, we obviously have
\beq
\label{p3214}
J_4:=\txinto g'(\br)\,z\,y_t\dx\ds\,\le\,\gamma\txinto y_t^2\dx\ds\,+\,\frac c\gamma\txinto z^2\dx\ds\,.
\eeq

Also, using \pier{\eqref{stabu2}} and the global bounds once more, we obtain that
\begin{align}
\label{p3215}
J_5:&=\txinto(\mula-\bmu)\,(g'(\rla)-g'(\br))\,y_t\dx\ds\nonumber\\[1mm]
&\le\,c\tint\|\mula(s)-\bmu(s)\|_6\,\|\rla(s)-\br(s)\|_3\,\|y_t(s)\|_2\,\ds\nonumber\\[1mm]
&\le\gamma\txinto y_t^2\dx\dx\,+\,\frac c\gamma\,\lambda^4\,.
\end{align}

Finally, invoking \eqref{taylor4} and Young's inequality, we have the estimate
\begin{align}
\label{p3216}
J_6:&=\,-\txinto \left(B[\rla]-B[\br]-\lambda DB[\br](\xi)\right)y_t\dx\ds\nonumber\\[1mm]
&\le\left\|B[\rla]-B[\br]-\lambda DB[\br](\xi)\right\|_{L^2(\gianni{Q_t})}\,\|y_t\|_{L^2(Q_t)}\nonumber\\[1mm]
&\le\,\gamma\txinto y_t^2\dx\ds\,+\,\gianni{\frac c\gamma \,\|y\|_{L^2(Q_t)}^2}
  +\frac c\gamma \,\left(\gianni R\left(\lambda\,\|h\|_{L^2(Q)}\right)\right)^2\,.
\end{align}

Thus, combining the estimates \eqref{p3210}--\eqref{p3216}, and choosing $\gamma=\frac 18$, we have
shown that, for every $t\in (0,T]$, we have the estimate
\begin{eqnarray}
&& \frac 14 \txinto \!|\yla_t|^2\dx\ds\,+\,\frac 12\,\|\yla(t)\|_{\gianni H}^2\,
\nonumber\\
&& \le\,c\left( \gianni{\tint\|y^\lambda(s)\|_H^2\ds}
+ \lambda^4
\,+\left(\gianni R\left(\lambda\,\|h\|_{L^2(Q)}\right)\right)^2\right)\,.
\label{p3217}
\end{eqnarray}

\vspace{4mm}
\underline{\sc Step 3:} \,\,\,We now add the estimates \eqref{p329} and \eqref{p3217}. It follows that,
with suitable global constants $c_1>0$ and $c_2>0$, we have for every $t\in (0,T]$ the estimate
\begin{align}
&\left\|\zla(t)\right\|_H^2\,+\left\|\zla\right\|^2_{\lzvt}\,+\left\|\yla(t)\right\|^2_H
\,+\left\|\yla_t\right\|^2_{\lzht}\nonumber\\[1mm]
&\le\,c_1\,Z(\lambda)\,+\,c_2\tint\left(1+\|\bmu_t(s)\|_V^2\right)\left(
\left\|\yla(s)\right\|_H^2\,+\left\|\zla(s)\right\|^2_H\right)\ds\,,
\end{align}
where we have defined, for $\lambda>0$, the function $Z$ by 
\beq
Z(\lambda):=\,\lambda^4\,+\,(\gianni R(\lambda\,\|h\|_{\lzq}))^2\,.
\eeq

Recalling that the mapping $\,s\mapsto \|\bmu_t(s)\|_V^2\,$ belongs to $L^1(0,T)$, we conclude from
Gronwall's lemma that, for every $t\in (0,T]$,
\begin{align}
&\left\|\yla\right\|_{H^1(0,t;H)}^2\,+\,\left\|\zla\right\|^2_{L^\infty(0,t;H)\cap \lzvt}\nonumber\\[1mm]
&\le\, c_1\,Z(\lambda) \exp\left(\gianni{c_2}\int_0^T\!\!\left(1+\|\bmu_t(s)\|_V^2\right)\ds\right)
\,\le\,c\,Z\left(\lambda\right)\,.
\end{align} 

Since $\,\lim_{\lambda\searrow 0}\,Z(\lambda)/\lambda^2=0\,$ (recall \eqref{taylor4}), we have
finally shown the validity of \eqref{Frechet}. This concludes the proof of the assertion.
\qed

\vspace{7mm}
We are now in the position to derive  the following result. 

\vspace{5mm}
{\sc Corollary 3.3:}  \quad{\em Let the general hypotheses} {\bf (A1)}--{\bf (A\gianni5)} {\em be
fulfilled and assume that $\bu\in\uad$ is a solution to the control problem} {\bf (CP)}
{\em with associated state $(\br,\bmu)=\cs (\bu)$. Then we have, for every $v\in \uad$,}
\begin{align}
\label{vug1}
&\beta_1\texinto (\br-\rho_Q)\,\xi\dx\dt\,+\,\beta_2\texinto(\bmu-\mu_Q)\,\eta\dx\dt
\nonumber\\[1mm]
&+\,\beta_3\texinto\bu\,(v-\bu)\dx\dt\,\ge\,0\,,
\end{align}
{\em where $(\xi,\eta)$ denotes the (unique) solution to the linearized system} 
\eqref{ls1}--\eqref{ls4} {\em for $h=v-\bu$.}

\vspace{3mm}
{\sc Proof:} \,\,\,Let $v\in\uad$ be arbitrary. Then $h=v-\bu$ is an admissible direction, since
$\bu+\lambda h\in\uad$ for $0< \lambda\le 1$. For any such $\lambda$, we have
\begin{align*}
0\, & \le\,\frac{J(\bu+\lambda h, \cs(\bu+\lambda h))-J(\bu,\cs(\bu))}\lambda
\\[1mm]
& \le\,\frac{J(\bu+\lambda h,\cs(\bu+\lambda h))-J(\bu,\cs(\bu+\lambda h))}\lambda
\\[1mm]
&\quad\,+\,\frac{J(\bu,\cs(\bu+\lambda h)-J(\bu,\cs(\bu))}\lambda\,\,.
\end{align*}

It follows immediately from the definition of the cost functional $J$ that the first summand on the right-hand
side of this inequality converges to $\,\,\texinto\beta_3\,\bu\,h\dx\dt\,\,$ as $\lambda\searrow 0$. For the
second summand, we obtain from Theorem 3.2 that
\begin{align*}
&\lim_{\lambda\searrow 0}\,\frac{J(\bu,\cs(\bu+\lambda h)-J(\bu,\cs(\bu))}\lambda
\\[1mm]
&=\,
\beta_1\texinto (\br-\rho_Q)\xi \dx\dt\, +\,\beta_2\texinto(\bmu-\mu_Q)\eta\dx\dt\,, 
\end{align*}

whence  the assertion follows.\qed

\section{Existence and first-order\\ necessary conditions of optimality}

\setcounter{equation}{0}
In this section, we derive the first-order necessary conditions of optimality for
problem {\bf (CP)}. We begin with an existence result.

\vspace{5mm}
{\sc Theorem 4.1:} \quad {\em Suppose that the conditions} {\bf (A1)}--{\bf (A\gianni5)} {\em are satisfied. Then
the problem} {\bf (CP)} {\em has a solution} $\bu\in\uad$.

\vspace{3mm}
{\sc Proof:} Let $\{u_n\}_{n\in\nz}\subset \uad$ be a minimizing sequence for {\bf (CP)}, and let $\{(\rho_n,
\mu_n)\}_{n\in\nz}$ be the sequence of the associated solutions to (\ref{ss1})--(\ref{ss4}). We then can infer from the global estimate \eqref{ssbounds2} the 
existence of a triple $(\bu,\br,\bmu)$ such that, for a suitable subsequence again indexed by~$n$,  
\begin{align*}
&u_n\to \bar u\,\quad\mbox{weakly star in }\,H^1(0,T;H)\cap L^\infty(Q),\\[1mm]
&\rho_n\to\bar\rho\,\quad\mbox{weakly star in }\,H^2(0,T;H)\cap
\pier{W^{1,\infty}(0,T;L^\infty(\oma))\cap H^1(0,T;V)},\\[1mm]
&\mu_n\to\bar\mu \,\quad\mbox{weakly star in }\,W^{1,\infty}(0,T;H)\cap H^1(0,T;V)\cap 
L^\infty(0,T;W).
\end{align*}
Clearly, we have that $\bu\in\uad$ and, by virtue of the Aubin-Lions lemma 
(cf.\ \cite[Thm.~5.1, p.~58]{Aubin}) and similar compactness results
(cf.\ \cite[Sect.~8, Cor.~4]{Simon}), 
\beq 
\rho_n\to \br\,\quad\mbox{strongly in }\,L^2(Q),
\label{eq:2.10}
\eeq
\gianni{whence also
$\rho_*\le \br\le\rho^*$ a.\,e.\ in $Q$ and
\begin{align*}
&B[\rho_n]\to B[\br] \quad\mbox{strongly in }\,\lzq,\\[1mm]
&\Phi(\rho_n)\to \Phi(\br) \quad\mbox{strongly in }\,L^2(Q),\quad \mbox{for }\,
\Phi\in\{F',g,g'\},
\end{align*} 
thanks to the general assumptions on $B$, $F$ and~$g$,
as well as the strong convergence
\beq
\label{eq:2.9}
\mu_n\to \bmu\,\quad\mbox{strongly in }\,\gianni{C^0([0,T];C^0(\overline\Omega))=C^0(\overline Q).}\\[1mm]
\eeq
From this, we easily deduce that
\begin{align*}
&g(\rho_n)\,\pt\mu_n\to g(\br)\,\pt\bmu \quad\mbox{weakly in }\,L^1(Q),\\[1mm]
&\mu_n\,g'(\rho_n)\,\pt\rho_n\to \bmu\,g'(\br)\,\pt\br\quad\mbox{weakly in }\,
L^1(Q)\,.
\end{align*}
}%
In summary, if we pass to the limit as $n\to\infty$ in the state equations 
(\ref{ss1})--(\ref{ss2}),
written for the triple $(u_n,\rho_n,\mu_n)$, we find that $(\br,\bmu)$ satisfies
\eqref{ss1} and \eqref{ss2}. Moreover, $\bmu\in L^\infty(0,T;W)$ satisfies the
boundary condition \eqref{ss3}, and it is easily seen that also the initial
conditions \eqref{ss4} hold true. In other words, we have $(\br,\bmu)=\cs(\bu)$,
that is, the triple $(\bu,\br,\bmu)$ is admissible for the control problem
{\bf (CP)}. From the weak sequential lower semicontinuity of the cost functional $J$
it finally follows that $\bar u$, together with $(\br,\bmu)=\cs(\bu)$, is
a solution to {\bf (CP)}. This concludes the proof. \qed

\vspace{5mm}
We now turn our interest to the derivation of first-order necessary optimality conditions for
problem {\bf (CP)}. To this end, 
we generally assume in the following that the hypotheses {\bf (A1)}--{\bf (A\gianni5)}
are satisfied and that $\bu\in\uad$ is an optimal control
with associated state
$(\br,\bmu)$, which has the properties (\ref{ssbounds1})--(\ref{ssbounds3}). We now aim to
eliminate $\xi$ \gianni{and $\eta$} from the variational inequality \eqref{vug1}. To this end, we employ
the adjoint state system associated with \eqref{ss1}--\eqref{ss4} for $\bu$, which is formally given by:
\begin{align}
  & -(1+2g(\br))\,p_t - g'(\br)\,\br_t\,p - \Delta p - g'(\br)\,q
  \nonumber
  \\[1mm]
  & \quad {} = \gian{\beta_2(\bmu-\mu_Q)}
  \quad\mbox{in }\,Q,
  \label{as1}
  \\[2mm]
  & -q_t + F''(\br)\,q - \bmu \,g''(\br)\,q + 
	g'(\br) \left( \bmu_t\, p-\bmu \,p_t \right) + DB[\br]^*(q)
  \nonumber
  \\[1mm]
  & \quad {} = \gian{\beta_1(\br-\rho_Q)}
  \quad\mbox{in }\,Q,
  \label{as2}
  \\[2mm]
  & \qquad \dn p=0\quad\mbox{on }\,\Sigma,
  \label{as3}
  \\[2mm]
  & \qquad p(T)=q(T)=0 \quad\mbox{in }\,\oma\,.
  \label{as4}
\end{align}
 
In \eqref{as2}, $DB[\br]^*\in {\cal L}(\lzq,\lzq)$ denotes the adjoint operator associated with the 
operator $\,DB[\br]\in {\cal L}(\lzq,\lzq)$, \gianni{thus} defined by the identity
\beq
\label{adjoint}
\texinto DB[\br]^*(v)\,w\dx\dt\,=\,\texinto v\,DB[\br](w)\dx\dt \quad\forall\,v,w\in \lzq\,.
\eeq
\gianni{As, for every $v\in L^2(Q)$, the restriction of $DB[\br](v)$ to $Q_t$ depends only on~$v|_{Q_t}$,
it follows that, for every $w\in L^2(Q)$,
the restriction of $DB[\br]^*(w)$ to~$Q^t \pier{{}=\Omega\times (t,T)}$ (see~\eqref{defQt}) depends only on~$w|_{Q^t}$.
Moreover, \eqref{B5} implies that
\beq
  \|DB[\br]^*(w)\|_{L^2(Q^t)} \leq C_B \|w\|_{L^2(Q^t)}
  \quad \forall\, w\in L^2(Q).
  \label{stimaDBstar}
\eeq
}%
We \gianni{also} note that in the case of the integral operator \eqref{intop2} it follows from Fubini's
theorem that $DB[\br]^*=DB[\br]=\gianni B$.

\vspace{2mm}
We have the following existence and uniqueness result for the adjoint system.

\vspace{5mm}
{\sc Theorem 4.2:} \quad 
{\em Suppose that} {\bf (A1)}--{\bf (A5)} 
{\em are fulfilled, and assume that $\bu\in\uad$ is a solution to the control problem} {\bf (CP)}
{\em with associated state $(\br,\bmu)=\cs (\bu)$.
Then the adjoint system} \accorpa{as1}{as4} {\em has a unique solution $(p,q)$ satisfying
\Beq
  p \in \H1H \cap \L\infty V \cap \L2W
  \aand
  q \in \H1H \,.
  \label{regsoluzadj}
\Eeq
}%

\vspace{3mm}
{\sc Proof:}  \,\,\,
Besides \juerg{of Young's} inequality, we make repeated use of  
the global estimates \accorpa{ssbounds1}{ssbounds2} and \eqref{ssbounds3} for $\br$ and~$\bmu$, 
without further reference.
Moreover, we denote by $c$ different positive constants 
that may depend on the given data of the state system and of the control problem; 
the meaning of $c$ may change between and even within lines.

\vspace{2mm}
We first prove uniqueness.
Thus, we replace the \rhs s of \eqref{as1} and \eqref{as2} by~$0$ and prove that $(p,q)=(0,0)$.
We add $p$ to both sides of \eqref{as1} and multiply by~$-p_t$.
At the same time, we multiply \eqref{as2} by~$q$.
Then we add the resulting equalities and integrate over~\pier{$Q^t=\Omega\times(t,T)$.}
As $g$ is nonnegative\juerg{, and thanks to \eqref{B5}, we~obtain that}
\Bsist
  && \bintQt |p_t|^2 \dx\ds
  \,+\, \frac 12 \, \normaV{p(t)}^2 
  \,+\, \frac 12 \iO |q(t)|^2 \dx
  \non
  \\
  && \leq \bintQt \bigl( - p - g'(\br)\, \br_t\, p - g'(\br)\, q \bigr)\, p_t \dx\ds
  \non
  \\
  && \quad {}
  + \bintQt \bigl(
    ( \bmu \,g''(\br) - F''(\br) )\, q 
    + \bmu \,g'(\br) \,p_t 
    - DB[\br](q) 
  \bigr) \,q \dx\ds
  \non
  \\
  && \quad {}
  - \bintQt g'(\br)\, \bmu_t\, p\, q \dx\ds
  \non
  \\
  && \leq \frac 12 \bintQt |p_t|^2 \dx\ds
  + c \bintQt ( p^2 + q^2 ) \dx\ds
  + c \bintQt |\bmu_t| \, |p| \, |q| \dx\ds \,.
  \non
\Esist
The last integral is estimated \juerg{as follows: employing the H\"older, Sobolev and Young inequalities, we have}
\Bsist
  && \bintQt |\bmu_t| \, |p| \, |q| \dx\ds 
  \leq \itt \norma{\bmu_t(s)}_3 \, \norma{p(s)}_6 \, \norma{q(s)}_2 \ds
  \non
  \\
  && \leq c \itt \bigl( \normaV{\bmu_t(s)}^2 \, \normaV{p(s)}^2 + \normaH{q(s)}^2 \bigr) \ds \,.
  \non
\Esist
As the function $s\mapsto\normaV{\bmu_t(s)}^2$ \juerg{belongs to} $L^1(0,T)$,
we can apply the backward version of Gronwall's lemma
to conclude that $(p,q)=(0,0)$.

\medskip

The existence of a solution to \accorpa{as1}{as4} is proved in several steps.

\step 1
We approximate $\br$ and $\bmu$ by functions $\rhoeps,\mueps\in C^\infty(\overline Q)$ satisfying
\accorpa{hpapprox}{conveps}
and look for a solution $(\peps,\qeps)$ to the following problem:
\Bsist
  && -(1+2g(\rhoeps))\,\peps_t - g'(\br)\,\br_t\,\peps - \Delta\peps - g'(\br)\,\qeps
  \non
  \\[1mm]
  && \quad {} = \gian{\beta_2(\bmu-\mu_Q)}
  \quad\mbox{in }\,Q,
  \label{as1e}
  \\[2mm]
  && -\qeps_t - \eps\,\Delta\qeps + F''(\br)\,\qeps - \bmu \,g''(\br)\,\qeps + 
  g'(\rhoeps) \left( \mueps_t\,\peps-\mueps\,\peps_t \right) 
  \non
  \\[1mm]
  && \quad {} + DB[\br]^*(\qeps)
  = \gian{\beta_1(\br-\rho_Q)}
  \quad\mbox{in }\,Q,
  \label{as2e}
  \\[2mm]
  && \qquad \dn\peps=\dn\qeps=0\quad\mbox{on }\,\Sigma,
  \label{as3e}
  \\[2mm]
  && \qquad \peps(T)=\qeps(T)=0 \quad\mbox{in }\,\oma\,.
  \label{as4e}
\Esist
We prove that \juerg{this} problem has a unique solution satisfying
\Beq
  \peps,\qeps \in \H1H \cap \L\infty V \cap \L2W \,.
  \label{regeps}
\Eeq
To this end, we present \accorpa{as1e}{as3e} as an abstract backward equation, 
\juerg{namely},
\Beq
  - \frac d{dt} \, (\peps,\qeps)(t)
  + \Aeps(t) \, (\peps,\qeps)(t)
  + (\Reps (\peps,\qeps))(t)
  = \feps(t),
  \label{backpbl}
\Eeq
in the framework of the Hilbert triplet $(\calV,\calH,\calV')$, where
\Beq
  \calV := V \times V
  \aand
  \calH := H \times H \,.
  \non
\Eeq
Notice that \eqref{backpbl}\juerg{, together with the regularity \eqref{regeps}, means that}
\Bsist
  && - \bigl( (\peps_t,\qeps_t)(t) , (v,w) \bigr)_{\calH}
  + \aeps\bigl( t;(\peps,\qeps)(t),(v,w) \bigr)
  \non
  \\
  && \quad {}
  + \bigl( (\Reps (\peps,\qeps))(t) , (v,w) \bigr)_{\calH}
  \,=\, \bigl( \feps(t) , (v,w) \bigr)_{\calH}
  \non
  \\[1mm]
  && \hbox{for every $(v,w)\in\calV$ and a.\,a.\ $t\in(0,T)$},
  \label{back}
\Esist
where $\aeps(t;\cdot,\cdot)$ is the bilinear form 
associated with the operator $\Aeps(t):\calV\to\calV'$; 
\juerg{moreover, $(\cdot,\cdot)_{\calH}$ denotes the inner product in~$\calH$
(equivalent to the usual one) that} one has chosen,
the embedding $\calH\subset\calV'$ being dependent on such a choice.
In fact, we will not use the standard inner product of~$\calH$,
which will lead to a nonstandard embedding $\calH\subset\calV'$.
We aim at applying first \cite[Thm.~3.2, p.~256]{Baiocchi},
in order to find a unique weak solution,
as we did for the linearized problem; then, we derive the full regularity required in~\eqref{regeps}.
We set, for convenience,
\Beq
  \phieps := \frac 1 {1+2g(\rhoeps)}
  \aand
  \psieps := \frac {\mueps \,g'(\rhoeps)} {1+2g(\rhoeps)}
  = \phieps \,\mueps\, g'(\rhoeps),
  \non
\Eeq
and choose a constant $\Meps$ such that
\Beq
  \phieps \leq \Meps , \quad
  |\psieps| \leq \Meps , \quad
  |\nabla\phieps| \leq \Meps,
  \aand
  |\nabla\psieps| \leq \Meps,
  \quad \hbox{a.\,e.\ in $Q$}.
  \non
\Eeq
Moreover, we introduce three parameters 
$\lambdaeps,\lambdaeps_1,\lambdaeps_2$,
whose values will be specified later on.
In order to transform our problem, we compute $\peps_t$ from \eqref{as1e}
and substitute in~\eqref{as2e}.
Moreover, we multiply \eqref{as1e} by~$\phieps$.
Finally, we add and subtract the same terms for convenience.
Then \accorpa{as1e}{as2e} is equivalent to the system
\Bsist
  && -\peps_t - \phieps\, \Delta\peps
  + \lambdaeps_1\, \peps
  \non
  \\
  && \quad {}
  - \lambdaeps_1 \,\peps
  - \phieps \,g'(\br)\,\br_t\, \peps - \phieps\, g'(\br)\, \qeps
  \,=\, \phieps\, \gian{\beta_2(\bmu-\mu_Q)},
  \non
  \\[2mm]
  && -\qeps_t
  - \eps\,\Delta\qeps
  + \psieps\, \Delta\peps
  + \lambdaeps_2 \,\qeps
  \non
  \\
  && \quad {}
  - \lambdaeps_2 \,\qeps
  + F''(\br)\,\qeps
  - \bmu \,g''(\br)\,\qeps
  + g'(\rhoeps)\, \mueps_t\,\peps
  \non
  \\
  && \quad {}
  + \psieps\, \bigl(    
    g'(\br)\,\br_t\, \peps
    + g'(\br)\, \qeps
    + \gian{\beta_2(\bmu-\mu_Q)}
  \bigr)
  \non
  \\
  && \quad {} + DB[\br]^*(\qeps)
  = \gian{\beta_1(\br-\rho_Q)} \,.
  \non
\Esist
By observing that
\Beq
  - \phieps \,\Delta\peps = - \div(\phieps\nabla\peps) + \nabla\phieps \cdot \nabla\peps\,,
  \non
\Eeq
and that the same holds true with $\psieps$ in place of~$\phieps$,
we see that the latter system, complemented with the boundary condition \eqref{as3e}, is equivalent~to
\Bsist
  && - \iO \peps_t(t)\, v \dx
  + \aeps_1(t;\peps(t),v)
  + \iO (\Reps_1(\peps,\qeps))(t)\, v \dx
  \non
  \\
  && \quad {}
  = \iO \pier{\phieps}(t) \, \gian{\beta_2(\bmu-\mu_Q)(t)} \, v \dx
  \non
  \\
  && - \iO \qeps_t(t)\, w \dx
  + \aeps_2(t;\pier{(\peps(t),\qeps(t))}, w)
  + \iO (\Reps_2(\peps,\qeps))(t) \,w \dx
  \non
  \\
  && \quad {}
  = \pier{- \iO \pier{\psieps}(t) \, \gian{\beta_2(\bmu-\mu_Q)(t)} \, w \dx} +\iO \gian{\beta_1(\br-\rho_Q)(t)} \, w \dx
  \non
\Esist
for every $(v,w)\in\calV$ and a.\,a.\ $t\in(0,T)$,
where the forms $\aeps_i$ are defined below
and the operators $\Reps_i$ account for all the other terms on the \lhs s of the equations.
We~set, for every $t\in[0,T]$ and $\hat v,\hat w,v,w\in V$,
\Bsist
  && \aeps_1(t;\hat v,v)
  := \iO \bigl(
    \phieps(t) \,\nabla\hat v \cdot \nabla v
    + (\nabla\phieps(t) \cdot \nabla\hat v)\, v
    + \lambdaeps_1 \,\hat v\, v
  \bigr) \dx\,,
  \non
  \\
  && \aeps_2(t;(\hat v,\hat w),w)
  \non
  \\
  && \quad {}
  := \iO \bigl(
    \eps\, \nabla\hat w \cdot \nabla w
    - \psieps(t) \,\nabla\hat v \cdot \nabla w
    - ( \nabla\psieps(t) \cdot \nabla\hat v)\, w
    + \lambdaeps_2\, \hat w \,w
  \bigr) \dx \,.
  \non
\Esist
Now, we choose the values of $\lambdaeps_i$ and of the further parameter $\lambdaeps$
\juerg{in such a way as} to guarantee some coerciveness.
\juerg{Putting} $\alpha:=1/(1+2\sup g)$, we have that
\Bsist
  && \aeps_1(t;v,v)
  \geq \iO \bigl(
    \alpha \,|\nabla v|^2 - \Meps\, |\nabla v| \, |v| + \lambdaeps_1\, v^2
  \bigr) \dx
  \non
  \\
  && \geq \iO \bigl(
    \alpha\, |\nabla v|^2
    - {\smfrac\alpha 2} \, |\nabla v|^2 
    - {\smfrac {\Meps^2}{2\alpha}} \, v^2 
    + \lambdaeps_1\, v^2
  \bigr) \dx \,.
  \non
\Esist
Therefore, the choice $\lambdaeps_1:= \frac\alpha 2 + \frac {\Meps^2}{2\alpha}$ yields
\Beq
  \aeps_1(t;v,v) \,\geq \,\frac\alpha 2 \, \normaV v^2
  \quad \hbox{for every $v\in V$ and $t\in[0,T]$.}
  \non  
\Eeq
Next, we deal with $\aeps_2$\juerg{. We have, for every $v,w\in V$ and $t\in[0,T]$,}
\Bsist
  && \hskip -1pt
  \aeps_2(t;(v,w),w) 
  \geq \iO \bigl(
    \eps |\nabla w|^2
    - \Meps |\nabla v| \, |\nabla w|
    - \Meps |\nabla v| \, |w|
    + \lambdaeps_2 w^2
  \bigr)\juerg{\dx}
  \non  
  \\
  && \geq \iO \bigl(
    \eps |\nabla w|^2
    - {\smfrac\eps 2 |\nabla w|^2}
    - {\smfrac {\Meps^2}{2\eps}} \, |\nabla v|^2
    - {\smfrac {\Meps^2}{2\eps}} \, |\nabla v|^2
    - {\smfrac\eps 2} \, |w|^2
    + \lambdaeps_2 w^2
  \bigr)\juerg{\dx}
  \non    
  \\
  && = \iO \bigl(
    {\smfrac\eps 2 |\nabla w|^2}
    + (\lambdaeps_2 - {\smfrac\eps 2}) w^2
    - {\smfrac {\Meps^2}\eps} \, |\nabla v|^2
  \bigr)\juerg{\dx\,,}
  \non  
\Esist
and the choice $\lambdaeps_2:=\eps$ leads to
\Beq
  \aeps_2(t;(v,w),w) 
  \,\geq\, \frac\eps 2 \, \normaV w^2 - \frac {\Meps^2}\eps \, \normaV v^2 \,.
  \non
\Eeq
Therefore, if we choose $\lambdaeps$ such that
$\lambdaeps \, \frac\alpha 2 - \frac {\Meps^2}\eps \geq \frac \eps 2$, \juerg{then}
we obtain
\Beq
  \lambdaeps \aeps_1(t;v,v)
  + \aeps_2(t;(v,w),w)
  \geq \frac\eps 2 \bigl( \normaV v^2 + \normaV w^2 \bigr)
  \non
\Eeq
for every $(v,w)\in\calV$ and $t\in[0,T]$.
Hence, if we define $\aeps:[0,T]\times\calV\times\calV\to\erre$ by setting
\Beq
  \aeps(t;(\hat v,\hat w),(v,w))
  := \lambdaeps_1 \aeps_1(\hat v,v)
  + \aeps_2(t;(\hat v,\hat w),w)\,,
  \non
\Eeq
\juerg{then} we obtain a time-dependent continuous bilinear form that is coercive on~$\calV$ 
(endowed whith its standard norm), 
uniformly with respect to $t\in[0,T]$.
Moreover, $\aeps$~depends smoothly on~$t$, and \accorpa{as1e}{as3e} is equivalent~to
\Bsist
  && - \iO \bigl( 
    \lambdaeps \peps_t(t) \, v
    + \qeps_t(t) \, w
  \bigr)\juerg{\dx}
  + \aeps \bigl( t;(\peps(t),\qeps(t)),(v,w) \bigr)
  \non
  \\
  && \quad {}
  + \iO \bigl\{
    \lambdaeps \,(\Reps_1(\peps,\qeps))(t) \, v
    + (\Reps_2(\peps,\qeps))(t) \, w
  \bigr\}\juerg{\dx}
  \non
  \\
  && = \iO \bigl(
    (\lambdaeps\,\phieps - \psieps)\juerg{(t)}\, \gian{\beta_2(\bmu-\mu_Q)(t)} \, v
    + \gian{\beta_1\,(\br-\rho_Q)(t)} \, w
  \bigr)\juerg{\dx}
  \non
\Esist
for every $(v,w)\in\calV$ and a.\,a.\ $t\in(0,T)$.
Therefore, the desired form \eqref{back} is achieved
if we choose the scalar product in~$\calH$ as follows:
\Beq
  \bigl( (\hat v,\hat w) , (v,w) \bigr)_{\calH}
  := \iO (\lambdaeps\,\hat v\, v + \hat w\, w) \dx
  \quad \hbox{for every $(\hat v,\hat w),(v,w)\in\calH$}.
  \non
\Eeq
Notice that this leads to the following nonstandard embedding $\calH\subset\calV'$:
\Beq
  {}_{\calV'}\< (\hat v,\hat w) , (v,w) >_{\calV}
  = \bigl( (\hat v,\hat w) , (v,w) \bigr)_{\calH}
  = \lambdaeps {}_{V'}\<\hat v , v >_V + {}_{V'}\<\hat w , w >_V
  \non
\Eeq
for every $(\hat v,\hat w)\in\calH$ and $(v,w)\in\calV$,
provided that the embedding $H\subset V'$ is the usual one, i.\,e.,
\juerg{corresponds to} the standard inner product of~$H$.
As the remainder, given by the terms $\Reps_1$ and $\Reps_2$,
satisfies the backward analogue of \accorpa{baioA}{baioB}
(see also \eqref{stimaDBstar}),
the quoted result of \cite{Baiocchi} can be applied,
and problem \accorpa{as1e}{as4e} has a unique solution satisfying
\Beq
  (\peps,\qeps) \in \H1{\calV'} \cap \L2\calV .
  \non
\Eeq
Moreover, if we move the remainder of \eqref{backpbl} to the \rhs,
we see that 
$$\pier{{}-\frac d{dt}(\peps,\qeps)+\Aeps(\peps,\qeps)\in\L2\calH .}$$ 
Therefore, by also accounting for~\eqref{as4e},
we deduce that $(\peps,\qeps)\in\H1\calH$
as well as $\Aeps(\peps,\qeps)\in\L2\calH$.
Hence, \pier{we have that} $\peps,\qeps\in\L2W$, by standard elliptic regularity.

\step 2
We add $\peps$ to both sides of \eqref{as1e} and multiply by~$-\peps_t$.
At the same time, we multiply \eqref{as2e} by~$\qeps$.
Then, we sum up and integrate over~$Q^t$. 
As $g\geq0$, we easily obtain that
\Bsist
  && \frac 12 \, \normaV{\peps(t)}^2
  + \bintQt |\peps_t|^2 \dx\ds
  + \frac 12 \iO |\qeps(t)|^2 \dx
  + \eps \bintQt |\nabla\qeps|^2  \dx\ds
  \non
  \\
  && \leq c \bintQt |\peps| \ |\peps_t| \dx\ds
  + c \bintQt |\qeps| \ |\peps_t| \dx\ds
  + c \bintQt |\qeps|^2 \dx\ds
  \non
  \\
  && \quad {}
  + c \bintQt |\mueps_t| \, |\peps| \, |\qeps| \dx\ds
  + \bintQt |DB[\br]^*(\qeps)| \, |\qeps| \dx\ds
 \pier{ {}+ c \,\norma\peps_{L^2(Q^t)}^2
  + c }\,.
  \non
\Esist
Just two \juerg{of the} terms on the \rhs\ need some \pier{treatment}.
We have
\Bsist
  && \bintQt |\mueps_t| \, |\peps| \, |\qeps| \dx\ds
  \leq \itt \norma{\mueps_t(s)}_3 \, \norma{\peps(s)}_6 \, \norma{\qeps(s)}_2 \ds
  \non
  \\
  && \leq c \itt \normaV{\peps(s)}^2 \ds
  + c \itt \norma{\mueps_t(s)}_3^2 \norma{\qeps(s)}_2^2 \ds\,,
  \non
\Esist
and we observe that the function $s\mapsto\norma{\mueps_t(s)}_3^2$
\juerg{belongs to}   $L^1(0,T)$, by \eqref{hpapprox}.
Moreover, the Schwarz inequality and \eqref{stimaDBstar} immediately yield that
\Beq
  \bintQt |DB[\br]^*(\qeps)| \, |\qeps| \dx\ds
  \leq C_B \,\norma\qeps_{L^2(Q^t)}^2 \,.
  \non
\Eeq
Therefore, we can apply the backward version of Gronwall's lemma to obtain that
\Beq
  \norma\peps_{\H1H\cap\L\infty V}
  + \norma\qeps_{\H1H}
  + \eps^{1/2} \norma\qeps_{\L2V}
  \leq c \,.
  \label{primastimaadj}
\Eeq
By comparison in \eqref{as1e}, we see that $\Delta\peps$ is bounded in~$\LQ2$.
Hence,
\Beq
  \norma\peps_{\L2W} \leq c \,.
  \label{daprimastimaadj}
\Eeq

\step 3
We multiply \eqref{as2e} by $-\qeps_t$ and integrate over~$Q^t$.
We obtain
\Bsist
  && \bintQt |\qeps_t|^2 \dx\ds \pier{{}+ \frac \eps 2 \int_\Omega |\nabla \qeps(t)|^2 \dx }
  \non
  \\
  && \leq c \bintQt |\qeps| \, |\qeps_t| \dx\ds
  + c \itt \norma{\mueps_t(s)}_3 \, \norma{\peps(s)}_6 \, \norma{\qeps_t(s)}_2 \ds
  \non
  \\
  && \quad {}
  + c \bintQt |\peps_t| \, |\qeps_t| \dx\ds
  + \bintQt |DB[\br]^*(\qeps)| \, |\qeps_t| \dx\ds \,.
  \non
\Esist
\pier{Thanks to} \eqref{stimaDBstar} once more, we deduce that
\Bsist
  && \frac 12 \bintQt |\qeps_t|^2 \dx\ds  \pier{{}+ \frac \eps 2 \int_\Omega |\nabla \qeps(t)|^2 \dx }
  \non
  \\
  && \leq c \bintQt |\qeps|^2 \dx\ds
  + c \itt \norma{\mueps_t(s)}_3^2 \, \normaV{\peps(s)}^2 \ds
  \non
  \\
  && \quad {}
  + c \bintQt |\peps_t|^2 \dx\ds \,.
  \non
\Esist
Thus, \eqref{hpapprox} and \eqref{primastimaadj} imply that
\Beq
\label{secondastimaadj}
\norma{\qeps_t}_{\L2H}  \pier{{}+ \eps^{1/2} \norma\qeps_{\L\infty V} }  \leq c \,.
\Eeq

\step 4
Now, we let $\eps$ tend to zero and construct a solution to~\accorpa{as1}{as4}.
By \accorpa{primastimaadj}{secondastimaadj} we have, at least for a subsequence,
\Bsist
  & \peps \to p
  & \quad \hbox{weakly star in $\H1H\cap\L\infty V\cap\L2W$}\,,
  \non
  \\
  & \qeps \to q
  & \quad \hbox{weakly in $\H1H$}\,,
  \non
  \\
  & \pier{ \eps \qeps \to 0}
  & \quad \pier{\hbox{strongly in $\L\infty V$}\,,}
  \non
\Esist
for some pair $(p,q)$ satisfying the regularity requirements \eqref{regsoluzadj}.
By accounting for \eqref{conveps} and the Lipschitz continuity of $g$ and~$g'$,
it is \sfw\ to see that $(p,q)$ is a solution to problem \accorpa{as1}{as4}.
This completes the proof.\qed

\vspace{5mm}
{\sc Corollary 4.3:} {\em Suppose that} {\bf (A1)}--{\bf (A\gianni5)} {\em are fulfilled, and assume that
$\bu\in\uad$ is an optimal control of {\bf (CP)} with associated state $(\br,\bmu)=\cs(\bu)$ and
adjoint state $(p,q)$. Then it holds the variational inequality}
\beq
\label{vug2}
\texinto (p+\beta_3\,\bu)(v-\bu)\dx\dt\,\ge\,0\quad\forall \,v\in\uad\,.
\eeq

\vspace{3mm}
{\sc Proof:}  \,\,\,
We fix $v\in\uad$ and choose $h=v-\bu$.
Then, we write the linearized system \accorpa{ls1}{ls4}
and multiply the equations \eqref{ls1} and \eqref{ls2} by $p$ and~$q$, respectively.
At the same time, we consider the adjoint system
and multiply the equations \eqref{as1} and \eqref{as2} by $-\eta$ and~$-\xi$, respectively.
Then, we \juerg{add all the equalities obtained in this way  and} integrate over~$Q$.
Many terms cancel out.
In particular, this happens for the contributions given by the Laplace operators, 
due to the boundary conditions \eqref{ls3} and~\eqref{as3},
\juerg{as well as} for the terms involving $DB[\br]$ and $DB[\br]^*$, by the definition of adjoint operator
(see \eqref{adjoint}).
Thus, it remains 
\Bsist
  && \intQ \bigl(
    2g'(\br)\, \br_t\, \eta\, p 
    + (1+2g(\br)) \,\eta_t\, p
    + (1+2g(\br)) \,\eta\, p_t
  \bigr) \dx\dt
  \non
  \\
  && \quad {}
  + \intQ \bigl(
    \bmu_t \,g'(\br)\, \xi\, p
    + \bmu \,g''(\br)\, \br_t \,\xi\, p
    + \bmu \,g'(\br)\, \xi_t\, p
    + \bmu\, g'(\br)\, \xi\, p_t
  \bigr) \dx\dt
  \non
  \\
  && \quad {}
  + \intQ (\xi_t \,q + \xi\, q_t) \dx\dt
  \non
  \\
  && = \intQ \bigl(
    (v-\bu)\, p 
    - \gian{\beta_2 (\bmu - \mu_Q)} \, \eta
    - \gian{\beta_1 (\br - \rho_Q)} \, \xi
  \bigr) \dx\dt
  \non\Esist
Now, we observe that the \juerg{expression on the} \lhs\ coincides with
\Beq
  \intQ \pt \bigl\{
    (1+2g(\br))\, \eta\, p
    + \bmu\, g'(\br)\, \xi\, p
    + \xi q
  \bigr\} \dx\dt \,.
  \non
\Eeq
Thus, it vanishes, due to the initial and final conditions \eqref{ls4} and~\eqref{as4}.
This implies that
\Beq
  \intQ \bigl(
    \beta_1 (\br - \rho_Q)\, \eta
    + \beta_2 (\bmu - \mu_Q)\, \xi
  \bigr) \dx\dt
  = \intQ (v-\bu)\, p \dx\dt \,.
  \non
\Eeq
Therefore, \eqref{vug2} follows from~\eqref{vug1}.\qed

\vspace{5mm}
{\sc Remark 4:} \,\,\,The variational inequality \eqref{vug2} forms together with the state system
\eqref{ss1}--\eqref{ss4} and the adjoint system \eqref{as1}--\eqref{as4} the system of first-order necessary
optimality conditions for the control problem {\bf (CP)}. Notice that in the case $\beta_3>0$ the function
$\,-\beta_3^{-1}p\,$ is nothing but the $L^2(Q)$ orthogonal projection of $\bu$ onto $\uad$.

 
\end{document}